\documentclass[10pt,a4paper,reqno]{amsart}
\usepackage{amsfonts}
\usepackage{amsthm}
\usepackage{amsmath}
\usepackage{amscd}
\usepackage[utf8]{inputenc}
\usepackage{t1enc}
\usepackage[mathscr]{eucal}
\usepackage{indentfirst}
\usepackage{graphicx}
\usepackage{graphics}
\usepackage{pict2e}
\usepackage{epic}
\usepackage{enumitem}
\usepackage{float}
\usepackage{MnSymbol}
\usepackage{multirow}
\usepackage{shuffle} %the symbol shuffle product
\usepackage[table,xcdraw]{xcolor} %table background color
\usepackage{hyperref}
\numberwithin{equation}{section}
\usepackage[margin=2cm]{geometry}
\usepackage{epstopdf} 
\usepackage[noadjust]{cite}

\usepackage{tikz}
\usetikzlibrary{decorations.markings,calc}
\def\centerarc[#1](#2)(#3:#4:#5)% Syntax: [draw options] (center) (initial angle:final angle:radius)
{ \draw[#1] ($(#2)+({#5*cos(#3)},{#5*sin(#3)})$) arc (#3:#4:#5); }

\def\centerarcpath(#1)(#2:#3:#4)% Syntax: [draw options] (center) (initial angle:final angle:radius)
{ ($(#1)+({#4*cos(#2)},{#4*sin(#2)})$) arc (#2:#3:#4) }

\tikzset{middlearrow/.style={
		decoration={markings,
			mark= at position 0.5 with {\arrow{#1}} ,
		},
		postaction={decorate}
	}
}

\usepackage[utf8]{inputenc}
\usepackage[T1]{fontenc}

\usepackage{mmacells}

\mmaDefineMathReplacement[≤]{<=}{\leq}
\mmaDefineMathReplacement[≥]{>=}{\geq}
\mmaDefineMathReplacement[≠]{!=}{\neq}
\mmaDefineMathReplacement[→]{->}{\to}[2]
\mmaDefineMathReplacement[⧴]{:>}{:\hspace{-.2em}\to}[2]
\mmaDefineMathReplacement{∉}{\notin}
\mmaDefineMathReplacement{∞}{\infty}

\mmaSet{
	morefv={gobble=2},
	linklocaluri=mma/symbol/definition:#1,
	morecellgraphics={yoffset=1.9ex}
}

\theoremstyle{plain}
\newtheorem{theorem}{Theorem}[section]
\newtheorem{lemma}[theorem]{Lemma}
\newtheorem{corollary}[theorem]{Corollary}
\newtheorem{proposition}[theorem]{Proposition}

\theoremstyle{definition}
\newtheorem{definition}[theorem]{Definition}
\newtheorem{conjecture}[theorem]{Conjecture}
\newtheorem{remark}[theorem]{Remark}
\newtheorem{problem}[theorem]{Problem}
\newtheorem{example}[theorem]{Example}

\newcommand{\mycomment}[1]{}

\newcommand{\CMZV}[2]{\textsf{CMZV}^{#1}_{#2}}
\newcommand{\MZV}[2]{\textsf{MZV}^{#1}_{#2}}
\DeclareMathOperator{\Li}{Li}

\begin{document}
	
	\title{Iterated integrals and multiple polylogarithm at algebraic arguments}
	
	\author[Kam Cheong Au]{Kam Cheong Au}
	
	\address{Rheinische Friedrich-Wilhelms-Universität Bonn \\ Mathematical Institute \\ 53115 Bonn, Germany} 
	
	\email{s6kmauuu@uni-bonn.de}
	

		\begin{abstract} By introducing a generalized notion of multiple zeta values associated with an arbitrary finite subset $S\subset \mathbb{P}^1(\mathbb{C})$ and studying their transformation properties under rational functions, we show that multiple polylogarithms evaluated at roots of unity (cyclotomic multiple zeta values, CMZVs) can be equivalently expressed in terms of iterated integrals involving certain non-root-of-unity arguments.
		
		We apply this theory to elucidate previously unknown 
			$\mathbb{Q}$-linear relations among CMZVs: they come from non-trivial 
		solutions of certain $S$-unit equations in the function field of $\mathbb{P}^1(\mathbb{C})$, thereby attaining the motivic dimension for low levels and weights. We introduce a datamine of CMZVs that appears to be the first rigorous compilation of this kind in the literature.
		
			In addition, we formulate several non-trivial Galois descent conjectures for multiple polylogarithms and present applications to certain Apéry-type infinite series.
	\end{abstract}
	
	\maketitle
	
	\section{Introduction}
	Various aspects of the \textit{multiple polylogarithm}
	$$\Li_{s_1,\cdots,s_k}(a_1,\cdots,a_k) := \sum_{n_1>\cdots>n_k\geq 1}\frac{a_1^{n_1}\cdots a_k^{n_k}}{n_1^{s_1} \cdots n_k^{s_k}},$$
	are closely connected to deep arithmetic questions, including functional equations and their special values \cite{charlton2021clean,GONCHAROV1995197,zagier1986hyperbolic,brown2014single}. In this article, we investigate a largely unexplored aspect of these special values in the case where the parameters $a_i$ are roots of unity.\par
	
	In this case, the multiple polylogarithm is known as a \textit{colored multiple zeta value} or \textit{cyclotomic multiple zeta value} (CMZV). More precisely, when the $a_i$ are $N$-th roots of unity, the $s_i$ are positive integers, and $(a_1, s_1) \neq (1,1)$, the number $\Li_{s_1,\cdots,s_k}(a_1,\cdots,a_k)\in \mathbb{C}$ is called a CMZV of \textit{weight} $n = s_1+\cdots+s_k$ and \textit{level} $N$. The special case $N=1$ is the well-known \textit{multiple zeta value}. Denote the $\mathbb{Q}$-span of CMZVs of weight $n$ and level $N$ by $\CMZV{N}{n}$. \par
	
	Like classical MZVs, CMZVs carry natural shuffle and stuffle algebra structures, a feature that has attracted considerable attention in the literature  \cite{bigotte2002lyndon, ihara2006derivation, ZhaoStandard, borwein2001special, panzer2017parity, racinet2002doubles, zhao2008multiple}.
	In this work, we present a new perspective on CMZVs: we show that one remains within the world of CMZVs even when certain parameters $a_i$ are no longer roots of unity. To illustrate this phenomenon, recall that CMZVs admit an interpretation in terms of iterated integrals over roots of unity:
	$$\CMZV{N}{n} = \text{Span}_\mathbb{Q}\left\{\int_{1>x_1>\cdots>x_n>0} \frac{dx_1}{x_1-c_1} \cdots \frac{dx_n}{x_n-c_n} \Big| c_i^N\in  \{0,1\}, c_1\neq 1, c_n\neq 0\right\}.$$
	When $N=5$, we will see that the space $\CMZV{5}{n}$ can equivalently be described as follows (see Example \ref{level5Ex}; $\mu=e^{2\pi i/5}$):
	\begin{equation}\label{aux_13}\text{Span}_\mathbb{Q}\left\{\int_{1>x_1>\cdots>x_n>0} \frac{dx_1}{x_1-c_1} \cdots \frac{dx_n}{x_n-c_n} \Big| c_i\in \{0,1,\frac{\sqrt{5}+3}{2},\frac{\sqrt{5}+1}{2},-\mu-\mu^2-\mu^3,1+\mu\}, c_1\neq 1, c_n\neq 0\right\}.\end{equation}
	
	Allowing non-roots of unity to appear in the iterated-integral description of CMZVs has several important advantages. In particular, it enables us to
	\begin{enumerate}
		\item express elements of $\CMZV{N}{n}$ in terms of $\Li_{s_1,\cdots,s_k}(a_1,\cdots,a_k) $ with $|a_i|<1$, thereby enabling rapid numerical evaluation, in contrast with the slow convergence of the original defining series;
		\item uncover certain mysterious $\mathbb{Q}$-linear relations between them that are predicted to exist but do not follow from evident relations;
		\item formulate conjectures reflecting deep Galois-descent phenomena for multiple polylogarithms;
		\item translate a wide range of classical problems into the well-developed framework of CMZVs—illustrated in this article by certain infinite sums.
	\end{enumerate}
	
	We shall return to points (2)–(4) after stating our main theorem and outlining the methodology. While we will not address point (1) in this article, it could have potential applications to the calculation of certain Feynman integrals
	 \cite{broadhurst2014multiple, henn2017evaluating, ablinger2019numerical,henn2013analytic, henn2014two}.
	
	\subsection{Methodology, main result and some consequences}
	The following is a corollary of our main result (Corollary \ref{main_corollary}). We first introduce the abbreviation $$\int_0^1 \omega(c_1)\cdots \omega(c_n) := \int_{1>x_1>\cdots>x_n>0} \frac{dx_1}{x_1-c_1} \cdots \frac{dx_n}{x_n-c_n}.$$
	\begin{corollary}\label{intro_corollary_1}
		Let $S = \{0,\infty,1,e^{2\pi i/N},\cdots,e^{2\pi i (N-1)/N}\}$, and let $R(x)$ be a rational function such that $R^{-1}(R(S)) = S$ and $\{0,1,\infty\} \subset R(S)$. Then the iterated integral
		$$\int_0^1 \omega(c_1)\cdots \omega(c_n) \in \CMZV{N}{n},$$ provided that $c_1\neq 1, c_n\neq 0, c_i\in R(S).$
	\end{corollary}
	The earlier claim (\ref{aux_13}) follows from the above result by choosing a suitable $R$. As a further consequence, consider the \textit{generalized polylogarithm} $\Li_{s_1,\cdots,s_k}(a_1)$, defined as $\Li_{s_1,\cdots,s_k}(a_1,\cdots,a_k)$ when $a_2=\cdots=a_k=1$. For many non-root-of-unity values $z$, these polylogarithms turn out to be CMZVs—often in highly non-obvious ways—as illustrated in Table~\ref{polylog_table}.
	
	\begin{table}[h]
		\begin{tabular}{|c|l|}
			\hline
			Level $N$ & \multicolumn{1}{c|}{$z$} \\ \hline
			4 &
			$\frac{1\pm i}{2},\, 1\pm i$ \\ \hline
			5 &
			$\frac{1-\sqrt{5}}{2},\,
			\frac{\sqrt{5}-1}{2},\,
			\frac{3-\sqrt{5}}{2},\,
			\frac{2-\mu+\mu^2-2\mu^3}{5}$ \\ \hline
			6 &
			$\frac{1}{4},\,
			\frac{1}{3},\,
			-\frac{1}{2},\,
			-\frac{1}{3},\,
			\frac{8}{9},\,
			\frac{1}{9},\,
			-\frac{1}{8},\,
			\frac{3+i\sqrt{3}}{6},\,
			-\frac{i}{\sqrt{3}},\,
			\frac{1-i\sqrt{3}}{4}$ \\ \hline
			7 &
			$\frac{1}{2}\csc\!\left(\frac{3\pi}{14}\right),\,
			\frac{1}{4}\csc^2\!\left(\frac{3\pi}{14}\right),\,
			\sin^2\!\left(\frac{\pi}{7}\right)\sec^2\!\left(\frac{\pi}{14}\right),\,
			\sin\!\left(\frac{\pi}{7}\right)\sec\!\left(\frac{3\pi}{14}\right)$ \\ \hline
			8 &
			$\frac{2-\sqrt{2}}{4},\,
			\frac{2-\sqrt{2}}{2},\,
			3-2\sqrt{2},\,
			-\frac{1}{\sqrt{2}},\,
			\frac{4-3\sqrt{2}}{8},\,
			12\sqrt{2}-16,\,
			(1-\sqrt{2})i$ \\ \hline
			10 &
			$\frac{1}{5},\,
			\sqrt{5}-2,\,
			\frac{\sqrt{5}+1}{4},\,
			9-4\sqrt{5},\,
			\frac{2}{\sqrt{5}},\,
			5\sqrt{5}+11,\,
			20\sqrt{5}-44,\,
			\frac{5-5\sqrt{5}}{16},\,
			\frac{9-4\sqrt{5}}{5}$ \\ \hline
			12 &
			$-\frac{1}{\sqrt{3}},\,
			1-\sqrt{3},\,
			21-12\sqrt{3},\,
			\frac{3-3\sqrt{3}}{8},\,
			\frac{12-7\sqrt{3}}{24},\,
			\frac{9-5\sqrt{3}}{18},\,
			27-15\sqrt{3},\,
			97-56\sqrt{3}$ \\ \hline
		\end{tabular}
			\caption{\small Some values of $z$ for which $\Li_{s_1,\cdots,s_k}(z)$ is a level-$N$ CMZV.}
		\label{polylog_table}
	\end{table}
	\vspace{-8pt}
	
	Our proof of Corollary \ref{intro_corollary_1} proceeds by introducing a new space $\MZV{S}{n}$, where $S$ is an arbitrary finite subset of $\mathbb{P}^1 :=  \mathbb{C} \cup \{\infty\}$. Roughly, it is the space spanned by iterated integrals of the form (see Section \ref{MZVS_section} for the full definition) $$\frac{dx_1}{x_1-c_1} \cdots \frac{dx_n}{x_n-c_n},\qquad c_i\in S,$$
	where the paths of integration range over all admissible paths in $\mathbb{P}^1$ whose endpoints lie in $S$.  \par 
	Our main result, from which the above corollary follows, is Theorem~\ref{CMZV_main_theorem}.
	\begin{theorem}[Main theorem]
		Let $N\geq 3$ and $S = \{0,\infty,1,\mu,\cdots,\mu^{N-1}\}$, where $\mu = e^{2\pi i / N}$. Then $\CMZV{N}{n} = \MZV{S}{n}$.
	\end{theorem}
	A key advantage of the space $\MZV{S}{n}$ is that it has useful transformation properties when the set $S$ is replaced by its image $R(S)$ under a rational function $R$ (Proposition~\ref{generalStrans}). These properties are completely invisible on the $\CMZV{N}{n}$ side and constitute a central strength of our approach. \par
	
	The theory underlying this result, developed throughout Section~2, is nevertheless non-trivial. The main difficulty lies in handling iterated integrals taken along various paths whose endpoints may coincide with singularities of the integrand. 
	
	To address this issue, we draw on an established framework developed by Racinet~\cite{racinet2002doubles} and later reformulated by Zhao~\cite{zhao2016multiple}. We shall consider certain group-like elements in a formal power series ring formed from iterated integrals. The Hopf algebra structure of the ring provides a convenient language for the regularization setup. Our arguments also require the machinery of tangential base points at both finite and infinite points of $\mathbb{P}^1$; see Proposition~\ref{mainregtheorem}.
	
	While the development of the theory and the proof of the main results (Sections~2–4) are entirely theoretical, explicit examples and applications necessarily involve a substantial amount of computation. These aspects are explored in Sections~5, 6, and~7.
	
	\subsection{Non-standard relations}
	A deep result of Deligne and Goncharov provides an upper bound for the motivic dimension of $\CMZV{N}{n}$: \begin{theorem}[\cite{deligne2010groupe, delignegroupes}]\label{deligne_bound}
		Let $D(n,N)$ be defined by
		$$1+\sum_{n=1}^\infty D(n,N) t^n = \begin{cases}(1-t^2-t^3)^{-1} \qquad &\text{ if } N = 1 \\ (1-t-t^2)^{-1} \qquad &\text{ if } N =2 \\ (1-at+bt^2)^{-1} \qquad &\text{ if } N \geq 3 \end{cases}$$
		where $a=\varphi(N)/2 + \nu(N), b=\nu(N)-1$. Here $\nu(N)$ denotes the number of distinct prime factors of $N$ and $\varphi$ is the Euler totient function. Then the motivic dimension of $\textsf{CMZV}^N_n$ is bounded above by $D(n,N)$. Moreover, if $N\in\{1,2,3,4,6,8\}$, then the motivic dimension is exactly $D(n,N)$.
	\end{theorem}
	
	Assuming Grothendieck's period conjecture \cite{grothendieck1966rham}, all $\mathbb{Q}$-linear relations between elements of $\CMZV{N}{n}$ should be motivic. While the dimension formula above does not provide an explicit description of such relations, it does specify how many independent relations must exist. \par
	
	For $N=1$ or $2$, it is widely believed that all $\mathbb{Q}$-relations should follow from stuffle and shuffle relations \cite{MZVdatamine}. For general $N$, the known mechanisms for producing $\mathbb{Q}$-relations among CMZVs were summarized by Zhao \cite{ZhaoStandard, zhao2008multiple, zhao2016multiple}. These include:
	\begin{itemize}
		\item shuffle relations, coming from the iterated-integral representation;
		\item stuffle relations, coming from the series definition; and
		\item distribution relations, reflecting symmetries among roots of unity.
	\end{itemize}
	In this work, we focus on identifying additional relations that are not explained by these standard mechanisms. Following Zhao, we refer to them as non-standard relations.
	
	For many composite levels $N$, numerical evidence indicates their existence. The smallest such example is $N=4$, where Zhao exploited the octahedral symmetry of the configuration $\{0,\pm 1,\pm i, \infty\} \subset \mathbb{P}^1$ to construct non-standard relations. For other values of $N$, however, they remain elusive. \par
	
	In Section~5, we introduce a new class of relations, which we call \textit{$S$-unit relations}, and show that they generate many new non-standard relations.
	
	\begin{definition}
		Let $S = \{0,\infty,1,\mu,\cdots,\mu^{N-1}\}$, where $\mu = e^{2\pi i /N}$. A rational function $R\colon\mathbb{P}^1\to \mathbb{P}^1$ is called \textit{$N$-unital} if $$R^{-1}(\{0,1,\infty\})\subset S.$$
	\end{definition}
	Finding such an $R$ is equivalent to solving an $S$-unit equation in the function field of $\mathbb{P}^1$: the condition says that both $R$ and $1-R$ have divisors supported on $S$. By a classical finiteness result \cite{silverman1984s,mason1983hyperelliptic}, only finitely many $N$-unital functions exist. \par
	
	We briefly describe the structure of an $S$-unit relation, postponing precise formulations to Section~5.2. Let $R_1,R_2,\cdots,R_r$ be $N$-unital functions such that
	\begin{itemize}
		\item $R_1(0) \in \{0,1,\infty\}$;
		\item $R_i(1) = R_{i+1}(0)$ for $i=1,\cdots,r-1$; and
		\item $R_r(1) = R_1(0)$.
	\end{itemize}
	Then the paths $R_i\circ [0,1]$ concatenate to form a loop in $\mathbb{P}^1$ based at a point of $\{0,1,\infty\}$. The $S$-unit relation is of the form
	$$\int_0^1 R_1^\ast w + \cdots + R_r^\ast w \equiv 0 \pmod{\text{product of lower weights}},$$
	where $w$ is a word in the alphabet $\{\frac{dx}{x},\frac{dx}{1-x}\}$. Note that the above integral is in $\CMZV{N}{}$ because the $R_i$ are $N$-unital. \par
	The author expects that the standard relations, together with $S$-unit relations, suffice to reach the motivic bound of Theorem~\ref{deligne_bound} for $N=6$ and $N=8$. The case $N=6$ is particularly noteworthy: although it has been studied by several authors, a complete set of relations does not seem to have been rigorously obtained before now. Ablinger \cite{ablinger2014iterated, ablinger2017discovering} developed \textit{ad hoc} methods to produce some, but not all, such relations. Independent numerical approaches motivated by applications in high-energy physics also revealed them \cite{henn2017evaluating, broadhurst2014multiple}. \par
	Our Corollary \ref{intro_corollary_1} also resolves questions concerning $\mathbb{Q}$-relations among certain generalized MZVs studied by Borwein and Broadhurst \cite{broadhurst2014multiple,borwein2001central,borwein2004experimentation,broadhurst2015multiple}. Two representative examples (discussed in Section~5.4) are
	$$\text{Multiple Deligne value}:\qquad \left\{\int_0^1 \omega(c_1)\cdots\omega(c_n)\bigg| c_i\in \{0,1,e^{2\pi i /6}\} \right\}$$
	$$\text{Multiple Landen value}:\qquad \left\{\int_0^1 \omega(c_1)\cdots\omega(c_n)\bigg| c_i\in \{0,1,\frac{1+\sqrt{5}}{2},\frac{3+\sqrt{5}}{2}\} \right\}$$
	
	\subsection{Datamine of CMZVs}
	The $S$-unit relations, combined with standard relations, allow us to rigorously construct many previously unknown $\mathbb{Q}$-linear relations, thereby reaching the motivic dimension for several levels and weights where this was previously inaccessible. Once the motivic dimension is attained, every CMZV can be expressed as a $\mathbb{Q}$-linear combination of a conjecturally independent set of basis constants.\par
	For $N=1,2$, a datamine already exists \cite{MZVdatamine}. For higher levels, however, no comparable systematic resource appears to be available in the literature\footnote{Partial or empirical databases by other authors do exist; see the end of Section~5.3 for an overview.}. We therefore introduce a Mathematica package\footnote{Available at \url{https://www.researchgate.net/publication/357601353}.}
	that provides explicit reductions for the following levels $N$ and weights $n$:
	\begin{itemize}
		\item $n\leq 5$ for $N=3$;
		\item $n\leq 6$ for $N=4$;
		\item $n\leq 4$ for $N=5$;
		\item $n\leq 5$ for $N=6$;
		\item $n\leq 4$ for $N=8$;
		\item $n\leq 3$ for $N=7,10,12$.
	\end{itemize}
	The package also makes the membership statement of Corollary~\ref{intro_corollary_1} explicit; the underlying algorithm is briefly described in Section~4. The datamine inspires several conjectures concerning Galois descent, which we describe next, and has been applied to a number of classical problems since the manuscript first appeared as a preprint.

	\subsection{Relations between polylogarithms and Galois descent}
	Table \ref{polylog_table} implies that $\Li_n(\alpha)$ is a CMZV of some level for many non-root-of-unity values $\alpha$. Once all linear relations among CMZVs of that level are known—i.e., once the motivic dimension is reached—one can rigorously verify identities among such polylogarithmic values. \par As an example, we give a conceptual (Section \ref{coxeterladder}) and effectively computable proof of Coxeter's famous ladder \cite{coxeter1935functions,lewin1991structural}:$$\text{Li}_2\left(\rho ^{20}\right)=2 \text{Li}_2\left(\rho ^{10}\right)+15 \text{Li}_2\left(\rho ^4\right)-10 \text{Li}_2\left(\rho ^2\right)+\frac{\pi ^2}{5} \qquad \rho = (\sqrt{5}-1)/2,$$
	revealing its close connection with the geometry of the regular icosahedron. Many other ladder identities admit a similar reinterpretation via CMZVs (see Section~6). Our goal here is not to establish new identities, but rather to demonstrate how CMZVs provide an interesting perspective for understanding them.\par
	
	Another direction suggested by Table~\ref{polylog_table} concerns Galois descent for multiple polylogarithms—a subtle and rather unexplored phenomenon. See \cite{charlton2026proof, xu2019evaluations, zheng2025proof, glanois2016motivic, charlton2024evaluation} and \cite[Conjecture~5]{broadhurst2014multiple} for examples and a discussion of Galois descent for CMZVs. Roughly speaking, it says that suitable Galois symmetrizations of CMZVs of a higher level appear to descend to a lower level. 
	
	As the simplest example, Table \ref{polylog_table} implies that $\Li_n(\frac{1\pm i}{2}) \in \CMZV{4}{n}$, and evidence in lower weights suggests that
		$$\Li_n\left(\frac{1+i}{2}\right) + \Li_n\left(\frac{1-i}{2}\right) \stackrel{?}{\in}  \CMZV{2}{n}.$$
	More generally, we make the following conjecture.
	\begin{conjecture}\label{galois_conj}For integers $n\geq 1$ and $k\geq 1$, we have the following. \\
		(a) $$\begin{aligned}\sum_{\substack{ s_1,\cdots,s_k\geq 1 \\ s_1+\cdots+s_k = n}} \left(\Li_{s_1,\cdots,s_k}\left(\frac{1+i}{2}\right) + \Li_{s_1,\cdots,s_k}\left(\frac{1-i}{2}\right)\right) &\stackrel{?}{\in}  \CMZV{2}{n}, \\
			\sum_{\substack{ s_1,\cdots,s_k\geq 1 \\ s_1+\cdots+s_k = n}} \left(\Li_{s_1,\cdots,s_k}(i) + \Li_{s_1,\cdots,s_k}(-i) \right)&\stackrel{?}{\in}  \CMZV{2}{n},\end{aligned}$$
		here the individual terms are in $\CMZV{4}{n}$. \\
		(b) We have $$\sum_{\substack{ s_1,\cdots,s_k\geq 1 \\ s_1+\cdots+s_k = n}} \left(\Li_{s_1,\cdots,s_k}(e^{\pi i /3}) + \Li_{s_1,\cdots,s_k}(e^{-\pi i /3})\right) \stackrel{?}{\in}  \CMZV{1}{n},$$
		here the individual terms are in $\CMZV{6}{n}$ (actually in $\CMZV{3}{n}$, see Example \ref{level6Ex_Deligne_value}).\\
		(c) Let $\rho = (\sqrt{5}-1)/2$. Then $$\Li_n(\rho^3) - \Li_n(-\rho^3) \stackrel{?}{\in} \CMZV{5}{n},$$
		here the individual terms are in $\CMZV{10}{n}$. \\
		(d) We have the following Galois descent from level $6$ to $2$ (see equation \eqref{level2conj} below): $$\Li_n\left(-\frac18\right) - \frac{3^{n-1}}{2^{n-2}}\Li_n\left(\frac14\right)  \stackrel{?}{\in} \CMZV{2}{n}.$$
	\end{conjecture}
	
	Part (b) of the conjecture has recently been checked up to weight $8$ by Sun and Zhou \cite{sun2026multiple}. Another Galois descent conjecture at level $4$ is the following.
	
	\begin{conjecture}\label{level4_conj2}
		For non-negative integers $a,b,c$, we have\footnote{Here we use notation that will be introduced in Section~2; in particular, $\shuffle$ denotes the shuffle product.}
		$$\int_0^1 \omega(0)^a \left[ \omega(i) + \omega(-i)\right] \left[\omega(1)^b \shuffle \omega(-1)^c\right]\in \CMZV{2}{1+a+b+c},$$
		equivalently, the value of the integral
		$$\int_0^1 \frac{x}{1+x^2} \log^a x \log^b (1-x) \log^c (1+x) dx \in \CMZV{2}{1+a+b+c}.$$
		Note that the left-hand sides are \textit{a priori} elements of $\CMZV{4}{}$.
	\end{conjecture}
	
	Here, we emphasize the importance of having a CMZV datamine: it enables us to formulate the above conjectures and to verify them in low weights.
	
	\subsection{Apéry-type series}
	The classical Apéry series, $$\sum _{n=1}^{\infty} \frac{(-1)^n}{n^3 \binom{2 n}{n}} = -\frac{2\zeta(3)}{5},$$ which played a role in Apéry’s proof of the irrationality of $\zeta(3)$, has inspired a vast literature devoted to discovering and rigorously proving similar identities. In particular, Sun \cite{sun2010conjectures, sun2021book} has conjectured numerous striking but highly non-trivial series identities, such as
	$$\begin{aligned}&\sum _{n=1}^{\infty } \frac{(-1)^{n-1} \left(10 H_n-\frac{3}{n}\right)}{n^3 \binom{2 n}{n}}=\frac{\pi ^4}{30}\\
		&\sum _{n=1}^{\infty } \frac{-102 H_n+3 H_{2 n}+\frac{28}{n}}{n^4 \binom{2 n}{n}}=-\frac{55}{18}\pi ^2 \zeta (3)\\
		&\sum _{n=1}^{\infty } \frac{H_{2 n}-H_n}{\binom{3 n}{n} \left(2^n n^2\right)} = -\frac{\pi  G}{2}+\frac{33 \zeta (3)}{32}+\frac{1}{24} \pi ^2 \log (2), \qquad G = \text{Catalan's constant}.\end{aligned}$$
	
	Many of these series, as well as numerous related examples, can be converted into CMZVs. When the corresponding weight and level are sufficiently small and a datamine is available (such as the one developed in this work), these identities can then be proved automatically. Although CMZV-based approaches to Apéry-type series are not new \cite{ablinger2014iterated, ablinger2017discovering, ablinger2019proving, ablinger2011harmonic, kalmykovBinomial, zhou2023sun, zhou2022hyper, au2020evaluation}, we illustrate how the theory developed in this article can unify this approach. 
	
	That said, the CMZV approach has inherent limitations. In particular, it cannot address more challenging cases in which the relevant weight or level is too large. For example, the following identities lie beyond the reach of current CMZV datamines:
	
	\vspace{7pt}
	\noindent\begin{minipage}{.5\linewidth}
		$$\begin{aligned}
			\sum _{n=1}^{\infty } \frac{2^n \left(-7 H_n+2 H_{2 n}+\frac{2}{n}\right)}{n^2 \binom{2 n}{n}}&=-\frac{\pi^2}{2} \log (2),\\
			\sum _{n=1}^{\infty } \frac{3^n \left(-8 H_n+6 H_{2 n}+\frac{5}{n}\right)}{n^2 \binom{2 n}{n}}&=\frac{26 \zeta (3)}{3},
		\end{aligned}$$
	\end{minipage}%
	\begin{minipage}{.5\linewidth}
		$$\begin{aligned}
			\sum _{n=1}^{\infty } \frac{6 H^{(2)}_{\left\lfloor n/2\right\rfloor }-\frac{(-1)^n}{n^2}}{n^2 \binom{2 n}{n}}&=\frac{13 \pi ^4}{1620}, \\
			\sum _{n=1}^{\infty } \frac{\binom{2 n}{n} \left(9 H_{2 n+1}+\frac{32}{2 n+1}\right)}{16^n (2 n+1)^2}&= 40 \beta (4)+\frac{5 \pi  \zeta (3)}{12}.
		\end{aligned}$$
	\end{minipage}
	\vspace{7pt}
	
	Such identities become more tractable when CMZV techniques are combined with hypergeometric transformation formulas arising from Wilf–Zeilberger pairs. Since this hybrid approach falls outside the scope of the present article, we refer the reader to \cite{chu2021further, chu2011dougall, au2025wilf, au2025multiple} for proofs and further discussion. Accordingly, the purpose of Section~7 is not to derive new results of this type, but rather to illustrate how far one can proceed using CMZV techniques alone.
	\vspace{4pt}
	
	This article is organized as follows. In Section~2, we establish notation and recall the necessary algebraic background on iterated integrals and tangential base points. Section~3 introduces the space $\MZV{S}{n}$, states our main theorem, and discusses some of its immediate consequences; the proof of this theorem is given in Section~4. Section~5 applies the developed framework to $S$-unit relations and introduces the datamine for CMZVs. The final two sections are devoted to applications: Section~6 concerns identities involving polylogarithms, while Section~7 discusses Apéry-type series.
	
	\section{Preliminaries}
	\subsection{Shuffle algebra}
	Let $X = \{x_1,\cdots,x_k\}$ be a finite set. Denote by $\mathbb{Q}\langle X\rangle$ the non-commutative polynomial ring over $\mathbb{Q}$ generated by $X$, and let $\mathbb{Q}\llangle X \rrangle$ be its completion (the ring of formal power series). Treating $X$ as an alphabet, let $X^*$ be the set of words (including the empty word) over $X$. \par
	The shuffle product $\shuffle$ on $\mathbb{Q}\langle X \rangle$ is defined inductively as follows:
	$$w \shuffle 1 = 1 \shuffle w = w, \qquad x w \shuffle y v = x(w\shuffle yv)+y(xu\shuffle v),$$
	for $w,v\in X^\ast$ and $x,y\in X$. One then distributes $\shuffle$ over addition and scalar multiplication; it is commutative and associative. Moreover, $\mathbb{Q}\langle X\rangle$ is a free algebra under the shuffle product, with the Lyndon words as a set of generators \cite[Theorem~6.1]{reutenauer1993free}. \par 
	\begin{lemma}\label{co-dim1-shuffle}
		Let $W$ be a codimension-one subspace of $X_1 := \text{Span}_\mathbb{Q}(x_1,\cdots,x_k) \subset \mathbb{Q}\langle X \rangle$. Let $V \subset \mathbb{Q}\langle X \rangle$ be a subspace such that
		\begin{enumerate}
			\item $V$ is closed under shuffle product,
			\item $V$ contains $X_1$, $1\in V$ and
			\item $w\mathbb{Q}\langle X \rangle \subset V$ for any $w\in W$.
		\end{enumerate} 
		Then $V=\mathbb{Q}\langle X \rangle$.
	\end{lemma}
	\begin{proof}
		Let $\{y_1,y_2,\cdots,y_{k-1}\}$ be a basis of $W$, and pick any $y_k\in X_1 -W$. If $Y = \{y_1,\cdots,y_{k-1},y_k\}$, then $\mathbb{Q}\langle X \rangle = \mathbb{Q}\langle Y \rangle$. Recall that the shuffle algebra $\mathbb{Q}\langle Y \rangle$ is freely generated by Lyndon words in $Y^\ast$, with the lexicographical order on $\{y_1,\cdots,y_{k-1},y_k\}$.\par By (3), the Lyndon words starting with $y_1,y_2,\cdots,y_{k-1}$ are in $V$. The only Lyndon word starting with $y_k$ is $y_k$ itself, and by (2), $y_k$ is also in $V$. Together, these words generate $\mathbb{Q}\langle Y \rangle$ under the shuffle product, so property (1) implies that $V$ is the whole space.
	\end{proof}
	
	The ring $\mathbb{C}\langle X\rangle$ and its completion $\mathcal{F}:=\mathbb{C}\llangle X\rrangle$ have a Hopf algebra structure. The coproduct is defined by $\Delta(x_i) = x_i\otimes 1 + 1\otimes x_i$, and the antipode $S$ is given by $S(x_1\cdots x_n) = (-1)^n x_n \cdots x_1$ for $x_i \in X$. \par
	For $x_1,x_2\in X$, we have the quotient maps $$\pi_1: \mathcal{F} \to \mathcal{F}/(x_1 \mathcal{F}), \qquad \pi_{1,2}: \mathcal{F}\to \mathcal{F}/(x_1 \mathcal{F} + \mathcal{F} x_2).$$ The Hopf algebra structure of $\mathcal{F}$ descends to these quotients. 
	
	\begin{proposition}\label{group_like_uniqueness}
		(a) Every group-like element in $\pi_{1}(\mathcal{F})$ is the image of a group-like element in $\mathcal{F}$ under $\pi_{1}$. Moreover, any two such elements $\Phi_1, \Phi_2$ are related by
		$$\Phi_2 = \exp(A x_1) \Phi_1$$ for some $A\in \mathbb{C}$.  \par
		(b) Every group-like element in $\pi_{1,2}(\mathcal{F})$ is the image of a group-like element in $\mathcal{F}$ under $\pi_{1,2}$. Moreover, any two such elements $\Phi_1, \Phi_2$ are related by
		$$\Phi_2 = \exp(A x_1) \Phi_1 \exp(B x_2)$$ for some $A,B\in \mathbb{C}$. 
	\end{proposition}
	\begin{proof}
		See \cite[Theorem~B.7]{zhao2016multiple}.
	\end{proof}
	
	Let $\mathbf{J}\in \mathcal{F}/(x_1 \mathcal{F} + \mathcal{F} x_2)$ be group-like, and let $\mathbf{I}$ be its unique group-like lift to $\mathcal{F}$ in which the coefficients of $x_1$ and $x_2$ are zero. One can recursively calculate the coefficients of $\mathbf{I}$ from those of $\mathbf{J}$ as follows \cite[Proposition~13.3.42]{zhao2016multiple}:
	
	Let $k,m,n\geq 0$ be integers and $\xi_i \in X$, with $\xi_1 \neq x_1$ and $\xi_k \neq x_2$, and set $\xi_1 \cdots \xi_k x_2^n = \xi_1 \cdots \xi_q$. 
	\begin{equation}\label{reg}\textbf{I}[x_1^m \xi_1 \cdots \xi_k x_2^n] = \begin{cases} 0 \qquad &\text{ if } mn=k=0, \\
			\textbf{J}[\xi_1 \cdots \xi_k]  \qquad &\text{ if } m=n=0, \\
			-\frac{1}{m}\sum_{i=1}^q \textbf{I}[x_1^{m-1} \xi_1 \cdots \xi_i x_1 \xi_{i+1}\cdots \xi_q]  \qquad &\text{ if } m>0, \\
			-\frac{1}{n}\sum_{i=1}^k \textbf{I}[\xi_1 \cdots \xi_{i-1} x_2 \xi_{i+1}\cdots \xi_k x_2^{n-1}]  \qquad &\text{ if } m=0,n>0,
	\end{cases}\end{equation}
	Throughout, we write $\textbf{I}[\omega]$ for the coefficient of $\omega$ in $\textbf{I}$.
	
	\subsection{Iterated integral}
	We briefly recall the required facts about iterated integrals \cite[Chap.~3]{gil2017multiple}, \cite{chen1977iterated}. For continuous functions $f_i(t)$ defined on $[a,b] \subset \mathbb{R}$, define inductively
	$$\int_a^b f_1(t) dt \cdots f_n(t) dt = \int_a^b f_1(u) du \cdots f_{n-1}(u) \int_a^u f_n(t) dt,$$
	When $n=1$, this is the usual definite integral $\int_a^b f_1(t)\,dt$; when $n=0$, we define its value to be $1$.
	This definition can be extended to a smooth manifold. Let $\gamma: [0,1]\to M$ be a path on a manifold $M$, and let $\omega_1,\cdots,\omega_n$ be differential $1$-forms on $M$. Then
	$$\int_\gamma \omega_1\cdots \omega_n := \int_0^1 f_1(t) dt \cdots f_r(t) dt$$
	where $\gamma^\ast \omega_i = f_i(t)\,dt$ is the pullback of $\omega_i$. If $f: N\to M$ is a differentiable map between two manifolds $N$ and $M$, then
	\begin{equation}\label{itintpullback}
		\int_{f\circ \gamma} \omega_1\cdots \omega_n = \int_\gamma f^\ast \omega_1 \cdots f^\ast\omega_n.
	\end{equation}
	Let $X = \{\omega_1,\omega_2,\cdots\}$ be a finite set of differential $1$-forms on a manifold $M$. By treating the elements of $X$ as letters, the iterated integral is a homomorphism under the shuffle product:
	$$\int_{\gamma} \omega_1\cdots \omega_n \int_{\gamma} \omega_{n+1}\cdots \omega_{n+m} = \int_\gamma \omega_1\cdots \omega_n \shuffle \omega_{n+1}\cdots \omega_{n+m}.$$
	
	Write $$\textbf{I}_\gamma(X) := \sum_{w\in X^\ast} \left(\int_\gamma w\right) w \in \mathbb{C}\llangle X\rrangle.$$ When the set $X$ is clear from the context, we write $\textbf{I}_\gamma$ for $\textbf{I}_\gamma(X)$.
	\begin{proposition}Let $X$ be a collection of continuous differential $1$-forms on a manifold $M$. \par
		(a) Let $\gamma$ be a path on $M$, $(\textbf{I}_\gamma)^{-1} = \textbf{I}_{\gamma^{-1}}$, with $\gamma^{-1}$ the reverse path of $\gamma$. \par
		(b) If $\gamma_1, \gamma_2$ are two paths on $M$ with $\gamma_1(1) = \gamma_2(0)$, then $\textbf{I}_{\gamma_2\gamma_1} = \textbf{I}_{\gamma_2}\textbf{I}_{\gamma_1}.$\par
		(c) $\textbf{I}_\gamma$ is a group-like element, i.e., $\Delta(\textbf{I}_\gamma) = \textbf{I}_\gamma\otimes \textbf{I}_\gamma$.\par
	\end{proposition}
	\begin{proof}
		The first two properties are easy to verify \cite[Theorem~3.19]{gil2017multiple}. For (c), see \cite[Prop.~13.3.13]{zhao2016multiple}. In fact, for any linear function $f:\mathbb{C}\llangle X \rrangle \to \mathbb{C}$, $\sum_{w\in X^\ast} f(w) w$ is group-like if and only if $f$ is a shuffle homomorphism. 
	\end{proof}
	
	For our applications to multiple polylogarithms, we will be mainly interested in differential $1$-forms on $\mathbb{C}$ of the shape \begin{equation}\label{aux_5}\omega(a):=\frac{dx}{x-a} \text{ for } a\in \mathbb{C}, \qquad \omega(\infty) := 0.\end{equation}
	If $\gamma: [0,1]\to \mathbb{C}$ is a path, then\footnote{Assuming that $\gamma$ does not pass through $a_1,\cdots,a_n$ except possibly at its endpoints.} the iterated integral $$\int_\gamma \omega(a_1)\cdots \omega(a_n)$$ 
	converges when $a_1\neq \gamma(1), a_n\neq \gamma(0)$. \par
	Recall our definition of \textit{multiple polylogarithm}: 
	$$\Li_{s_1,\cdots,s_k}(x_1,\cdots,x_k) =  \sum_{n_1>\cdots>n_k\geq 1}\frac{x_1^{n_1}\cdots x_k^{n_k}}{n_1^{s_1} \cdots n_k^{s_k}}, \qquad |x_1|<1,\quad |x_1x_2|<1,\quad\cdots,\quad |x_1\cdots x_k|<1.$$
	It can be rewritten as an iterated integral using words of the form $\omega(a)$. More precisely,
	\begin{equation}\label{polylogtoint}\Li_{s_1,\cdots,s_k}(x_1,\cdots,x_k) = (-1)^k \int_0^1 \omega(0)^{s_1-1}\omega(a_1)\cdots \omega(0)^{s_k-1}\omega(a_k) \qquad a_i = x_1^{-1}\cdots x_i^{-1}.\end{equation}
	
	\subsection{Tangential base point}\label{tangential_base_point_section}
	For the rest of this section, we use the following notation. 
	\begin{itemize}[leftmargin=*]
		\item Let $S = \{\infty,a_1,\cdots,a_k\}$ be a fixed finite subset of the extended complex plane $\mathbb{P}^1 := \mathbb{C} \cup \{\infty\}$.
		\item $X = \{\omega(a_1),\cdots,\omega(a_k)\}$ and the formal power series ring $\mathcal{F} = \mathbb{C}\llangle X \rrangle$.
		\item For $a\in \mathbb{P}^1$, $$\Omega(a) := \begin{cases}\omega(a) \quad &\text{ if }a\neq \infty, a\in S \\ -\omega(a_1)-\cdots-\omega(a_k) &\text{ if }a=\infty \\
			0  &\text{ if }a\notin S 
		\end{cases}.$$ 
		\item For $\textbf{I} \in \mathcal{F}$ and $u\in X^\ast$, we let $\textbf{I}[u]$ denote the coefficient of the monomial $u$ in $\textbf{I}$ and extend this notation linearly to $\mathcal{F}$.
	\end{itemize} 
	
	We first give two ``residue theorems'' for iterated integrals; their proofs are straightforward.
	\begin{lemma}[Residue at finite point]\label{radius0}
		Let $\rho(\varepsilon)$ be an arc of the circle centered at $a \in S$ with radius $\varepsilon$. Then, as $\varepsilon\to 0$\footnote{Throughout, the $O$-term for formal power series is interpreted coefficientwise.}, $$\textbf{I}_{\rho(\varepsilon)}(X) = \exp(A \omega(a)) + O(\varepsilon^{1/2}) = \exp(A\Omega(a))+ O(\varepsilon^{1/2})$$ for some $A\in \mathbb{C}  $. 
	\end{lemma}
	\begin{proof}
		After translating if necessary, we may assume that $a = 0$. Recall that $$\textbf{I}_{\rho(\varepsilon)}(X) = \sum_{\omega\in X^\ast} (\int_{\rho(\varepsilon)} \omega)\omega.$$ For any $\omega$ that contains a letter other than $\omega(0)$, the integral tends to $0$ with $\varepsilon$. Indeed, parameterize the arc $\rho(\varepsilon)$ by $\varepsilon e^{i\theta}$, where $\alpha \leq \theta \leq \beta$. Then $f^\ast \omega(0) = i\,d\theta$, while $f^\ast \omega(a) = \varepsilon e^{i\theta}/(\varepsilon e^{i\theta}-a)\,d\theta$ converges uniformly to $0$ if $a\neq 0$. Hence $$\int_{\rho(\varepsilon)}\omega(0)\cdots \omega(a) \cdots \omega(0) = \int_\alpha^\beta f^\ast\omega(0)\cdots f^\ast\omega(a) \cdots f^\ast\omega(0) = O(\varepsilon^{1/2}).$$
		Thus, only words of the form $\omega = \omega(0)^n$ remain. An explicit calculation gives $$\int_{\rho(\varepsilon)} \omega(0)^n = \frac{1}{n!} \left(\int_{\rho(\varepsilon)} \omega(0) \right)^n = \frac{1}{n!} \left(i(\beta-\alpha) \right)^n.$$ Therefore, the lemma holds with $A = i(\beta-\alpha)$.
	\end{proof}
	
	\begin{lemma}[Residue at $\infty$]\label{radiusinf}
		Let $\rho(R)$ be an arc of the circle centered at $a$ with radius $R$. Then, as $R\to \infty$, $$\textbf{I}_{\rho(R)}(X) = \exp(A (\omega(a_1)+\cdots+\omega(a_k))) + O(R^{-1})= \exp(-A\Omega(\infty)) + O(R^{-1/2})$$ for some $A\in \mathbb{C}$. 
	\end{lemma}
	\begin{proof}
		Parametrize the arc $\rho(R)$ by $Re^{i\theta}$, where $\alpha\leq\theta\leq \beta$. Then $f^\ast \omega(a) = iRe^{i\theta}/(Re^{i\theta}-a)\,d\theta$ converges uniformly to $i\,d\theta$ as $R\to \infty$, so  $$\int_{\rho(R)}\omega(c_1)\cdots \omega(c_n) = \int_\alpha^\beta f^\ast\omega(c_1)\cdots f^\ast\omega(c_n) \to \frac{(i(\beta-\alpha))^n}{n!} + O(R^{-1/2}).$$
		Therefore, the lemma holds with $A=i(\beta-\alpha)$.
	\end{proof}
	
	Let $\gamma: [0,1]\to \mathbb{P}^1$ be a piecewise smooth path that does not pass through $S$ except at its endpoints. We say that $\gamma$ is \textit{regular} at the endpoint $1$ if, for some $c\neq 0$, we have
	$$	\gamma(1-\varepsilon) =\begin{cases}
		\gamma(1) + c\varepsilon + O(\varepsilon^2), \quad &\gamma(1)\neq \infty \\ 
		\frac{c}{\varepsilon} + O(1) , \quad &\gamma(1)= \infty 
	\end{cases} \qquad\varepsilon \to 0.$$
	Similarly, we can define the notion of being regular at the endpoint $0$.
	\begin{theorem}
		Let $\gamma$ be a path from $\gamma(0)\notin S$ to a point $b = \gamma(1)\in S$ that is regular at $1$. Then 
		$$\lim_{\varepsilon \to 0} e^{-(\log \varepsilon) \Omega(b)} \textbf{I}_{\gamma|[0,1-\varepsilon]}(X)$$
		exists\footnote{That is, the limit exists for each fixed coefficient.} in $\mathbb{C}\llangle X \rrangle$. Here $\gamma|[0,1-\varepsilon]$ is the restriction of $\gamma$ to the interval $[0,1-\varepsilon]$.
	\end{theorem}
	\begin{proof}
		Abbreviate $\textbf{I}_{\gamma|[0,1-\varepsilon]}(X) = \textbf{I}_\varepsilon$. Let $\textbf{A}_\varepsilon = e^{-(\log \varepsilon) \Omega(b)} \textbf{I}_\varepsilon$. Note that $\textbf{A}_\varepsilon$ is group-like because both $e^{-(\log \varepsilon) \Omega(b)}$ and $\textbf{I}_\varepsilon$ are. Hence $\textbf{A}_\varepsilon[u\shuffle v] =\textbf{A}_\varepsilon[u]\textbf{A}_\varepsilon[v]$. Let \begin{equation}\label{aux_2}V = \{u\in \mathbb{Q}\langle X \rangle| \lim_{\varepsilon \to 0} \textbf{A}_\varepsilon[u] \text{ exists}\}.\end{equation}
		We need to show that $V=\mathbb{Q}\langle X \rangle$. Clearly, $V$ is a subspace closed under the shuffle product, and $1\in V$. Next, we show that $V$ contains all words of weight $1$.
		\begin{itemize}
			\item In the case $b\neq \infty$, assume that $b=a_1$. For $i\neq 1$, we have $\textbf{A}_\varepsilon [\omega(a_i)] = \int_{\gamma|[0,1-\varepsilon]} \frac{dx}{x-a_i}$, whose limit clearly exists because the endpoint $\gamma(1) = a_1\neq a_i$. For $i=1$, we have $$\textbf{A}_\varepsilon [\omega(a_1)] = -\log \varepsilon + \int_{\gamma|[0,1-\varepsilon]} \frac{dx}{x-a_1} = -\log \varepsilon + \log(\gamma(1-\varepsilon)-a_1) - \log(\gamma(0)-a),$$
			and since $\gamma(1-\varepsilon) = a_1 + c\varepsilon + O(\varepsilon^2)$ with $c\neq 0$ by the regularity assumption, the above limit exists. 
			\item The case $b=\infty$. Recall our convention $\Omega(\infty) = -\omega(a_1) - \cdots - \omega(a_k)$. We have $$\textbf{A}_\varepsilon [\omega(a_i)] = \log \varepsilon + \int_{\gamma|[0,1-\varepsilon]} \frac{dx}{x-a_i} =  \log \varepsilon + \log(\gamma(1-\varepsilon)-a_1) - \log(\gamma(0)-a),$$
			Since $\gamma(1-\varepsilon) = c\varepsilon^{-1} + O(1)$ with $c\neq 0$ by the regularity assumption, the limit exists. 
		\end{itemize}
		In both cases, we see that $V$ contains all words of weight $1$. If we can find a codimension-one subspace $W$ of the words of weight $1$ such that $w\mathbb{Q}\langle X\rangle\subset V$ for any $w\in W$, then Lemma \ref{co-dim1-shuffle} implies that $V = \mathbb{Q}\langle X \rangle$, which completes the proof. \par 
		In the case $b\neq \infty$, take $W = \text{Span}\{\omega(a_2),\cdots,\omega(a_k)\}$. For any $w\in W$ and $u\in w\mathbb{Q}\langle X\rangle$, we have $\textbf{A}_\varepsilon [u] = \textbf{I}_\varepsilon[u]$, and the limit exists because $u$ does not start with $\omega(a_1)$. \par
		In the case $b=\infty$, we can take $W = \text{Span}\{\omega(a_p)-\omega(a_q) | p\neq q\}$, which is the codimension-one subspace whose coefficients sum to zero. To see this, let $u = \omega(a_p)\omega(a_{i_2})\cdots \omega(a_{i_n}) := \omega(a_p) \theta$ and $v = \omega(a_q)\omega(a_{i_2})\cdots \omega(a_{i_n}) = \omega(a_q) \theta$. We compute $$\textbf{A}_\varepsilon[u] =  \sum_{j=1}^n \frac{(-\log \varepsilon)^j}{j!} \textbf{I}_\varepsilon[\omega(a_{i_{j+1}})\cdots \omega(a_{i_n})] + \textbf{I}_\varepsilon[\omega(a_p) \theta].$$
		Similarly, $$\textbf{A}_\varepsilon[v] =  \sum_{j=1}^n \frac{(-\log \varepsilon)^j}{j!} \textbf{I}_\varepsilon[\omega(a_{i_{j+1}})\cdots \omega(a_{i_n})] + \textbf{I}_\varepsilon[\omega(a_q) \theta].$$
		Subtracting, we obtain \begin{equation}\label{aux_4}\textbf{A}_\varepsilon[(\omega(a_p)-\omega(a_q))\theta] = \textbf{I}_\varepsilon[(\omega(a_p)-\omega(a_q))\theta].\end{equation}
		As $\omega(a_p)-\omega(a_q) = (\frac{1}{x-a_p} - \frac{1}{x-a_q}) dx$, the integrand is $O(1/x^2)$, so the limit of the right-hand side exists as $\varepsilon \to 0$. Thus, the claimed $W$ works. 
	\end{proof}
	
	\begin{remark}
		In the previous theorem, the case $b\neq \infty$ can also be proved by applying the projection that eliminates the ``divergent word'' from $\textbf{I}_{\gamma|[0,1-\varepsilon]}(X)$, taking the limit, and then appealing to Proposition \ref{group_like_uniqueness}; see \cite[Chap.~13]{zhao2016multiple} for details. However, it is not obvious how to adapt this approach to the case $b=\infty$. Our proof above applies to both cases.
	\end{remark}
	
	\begin{proposition}
		Let $\gamma$ be a path from $\gamma(0)\notin S$ to a point $b = \gamma(1)\in S$ that is regular at $\gamma(1)$. Then there exists a group-like element $\widetilde{\textbf{I}_\gamma} \in \mathbb{C}\llangle X \rrangle$ such that
		$$\lim_{\varepsilon \to 0}  \textbf{I}_{\gamma|[0,1-\varepsilon]}(X) = e^{(\log \varepsilon) \Omega(b)} \widetilde{\textbf{I}_\gamma} + O(\varepsilon^{1/2}).$$
	\end{proposition}
	\begin{proof}
		Let $\widetilde{\textbf{I}_\gamma}$ be the limit $\lim_{\varepsilon \to 0} e^{-(\log \varepsilon) \Omega(b)} \textbf{I}_{\gamma|[0,1-\varepsilon]}(X)$, which exists by the previous theorem. It is group-like because taking limits preserves this property. By general results on iterated integrals at singular points \cite[Lemma~3.3.20]{zhao2016multiple}, the error term for each degree-$n$ monomial in the above limit is in fact $O(\varepsilon (\log \varepsilon)^n)$, which is $O(\varepsilon^{1/2})$. \par 
		Alternatively, one can change $V$ in equation (\ref{aux_2}) to
		$$V = \{u\in \mathbb{Q}\langle X \rangle| L:=\lim_{\varepsilon \to 0} \textbf{A}_\varepsilon[u] \text{ exists and } \textbf{A}_\varepsilon[u] = L+O(\varepsilon^{1/2})\},$$
		and then argue as before.
	\end{proof}
	
	\begin{proposition}[Tangential base point]\label{mainregtheorem}
		Let $a=\gamma(0)$ and $b=\gamma(1)$ be the initial and terminal points of $\gamma$, respectively, and assume that $\gamma$ is regular at both endpoints. Then there exists a unique group-like element $\widetilde{\textbf{I}_\gamma} \in \mathbb{C}\llangle X \rrangle$ such that $\widetilde{\textbf{I}_\gamma}[\Omega(a)] = \widetilde{\textbf{I}_\gamma}[\Omega(b)] = 0$ and 
		$$\textbf{I}_{\gamma | (\varepsilon,1-\varepsilon)} = e^{(B+\log \varepsilon)\Omega(b)} \widetilde{\textbf{I}_\gamma} e^{(A-\log \varepsilon)\Omega(a)} + O(\varepsilon^{1/2}), \qquad \varepsilon\to 0,$$
		for some $A,B\in \mathbb{C}$. 
	\end{proposition}
	\begin{proof}
		Split $\gamma$ into the two paths $\gamma(t/2)$, $0\leq t\leq 1$, and $\gamma(1-t/2)$, $0\leq t\leq 1$, and then apply the above proposition to them. Uniqueness follows from Proposition \ref{group_like_uniqueness}.
	\end{proof}
	
	The bijective holomorphic maps on $\mathbb{P}^1$ are the Möbius transformations:
	$$R(x) = \frac{ax+b}{cx+d}, \qquad a,b,c,d\in \mathbb{C},\qquad ad-bc\neq 0.$$
	Recall that we defined $\omega(\infty) := 0$ in (\ref{aux_5}). With this convention, one checks that, for any Möbius transformation $R:\mathbb{P}^1 \to \mathbb{P}^1$ and any $a\in \mathbb{P}^1$, 
	$$R^\ast \omega(a) = \frac{R'}{R-a} dx = \omega(R^{-1}(a)) - \omega(R^{-1}(\infty)). $$
	For our fixed $S\subset \mathbb{P}^1$, let $G$ be its symmetry group, namely
	$$G:= \{ \text{Möbius transformation } R |  R(S) = S\}.$$
	The action of $G$ on $\mathbb{C}\llangle X \rrangle$ is defined by $g\omega(a) := (g^{-1})^\ast \omega(a) = \omega(g(a))-\omega(g(\infty))$ and extended to all of $\mathbb{C}\llangle X \rrangle$. There is another action\footnote{We distinguish this action from the one above by placing a dot in front.} of $G$ on the group-like elements of $\mathbb{C}\llangle X \rrangle$, as follows:
	$$g\cdot \sum_{w\in X^\ast} f(w) w := \sum_{w\in X^\ast} f(g^{-1}w) w,$$
	Here, the right-hand side is still a group-like element because $x\mapsto f(g^{-1}x)$ is a shuffle homomorphism. \par
	
	For a path $\gamma: [0,1]\to \mathbb{P}^1$, $\gamma$ is regular at both endpoints if and only if $g\circ \gamma$ is, because $g$ is a Möbius transformation. We consider the group-like elements $\widetilde{\textbf{I}_{g\circ \gamma}}$ and $\widetilde{\textbf{I}_\gamma} \in \mathbb{C}\llangle X \rrangle$ defined in Proposition \ref{mainregtheorem}.
	
	\begin{lemma}\label{group_trans_reg}
		Under the assumptions of the previous paragraph, there exist $A,B\in \mathbb{C}$ such that $$\widetilde{\textbf{I}_{g\circ \gamma}} = e^{B \Omega(g\gamma(1))} (g\cdot\widetilde{\textbf{I}_\gamma}) e^{A\Omega(g\gamma(0))}.$$
	\end{lemma}
	\begin{proof}
		If both $\gamma(0),\gamma(1)\notin S$, we need to show that $\widetilde{\textbf{I}_{g\circ \gamma}}  = g\cdot\widetilde{\textbf{I}_\gamma}$. In this case, all coefficients of $\widetilde{\textbf{I}_{g\circ \gamma}}$ are convergent integrals:
		\begin{align}\begin{split}\label{aux_3}\widetilde{\textbf{I}_{g\circ \gamma}} &= \sum_{w\in X^\ast} \left(\int_{g\circ \gamma} w\right)w = \sum_{w\in X^\ast} \left(\int_{\gamma} g^\ast w\right)w\\
				&= \sum_{w\in X^\ast} \left(\int_{\gamma} g^{-1} w\right)w = g\cdot \widetilde{\textbf{I}_\gamma},\end{split}\end{align}
		Next, consider the case $\gamma(0)\notin S$ and $\gamma(1)\in S$. We need to show that there exists $B\in \mathbb{C}$ such that $\widetilde{\textbf{I}_{g\circ \gamma}} = e^{B \Omega(g\gamma(1))} (g\cdot\widetilde{\textbf{I}_\gamma})$. Let $\textbf{A}$ and $\textbf{B}$ denote the left- and right-hand sides, respectively, and set
		$$V:= \{u\in \mathbb{Q}\langle X\rangle | \textbf{A}[u] = \textbf{B}[u]\},$$
		We need to show that $V = \mathbb{Q}\langle X\rangle$. Since $\textbf{A}$ and $\textbf{B}$ are group-like elements, $V$ is closed under the shuffle product and $1\in V$. \par
		In the case $g\gamma(1)\neq \infty$, say $g\gamma(1) = a_1$, let $W = \text{Span}\{\omega(a_2),\cdots,\omega(a_k)\}$. Then, for any $w\in W$ and $u\in w \mathbb{Q}\langle X \rangle$, $\widetilde{\textbf{I}_{g\circ \gamma}}[u]$ represents a convergent integral. Thus, (\ref{aux_3}) shows that $\widetilde{\textbf{I}_{g\circ \gamma}}[u] = (g\cdot \widetilde{\textbf{I}_\gamma})[u]$, so $u\in V$. \par
		In the case $g\gamma(1) = \infty$, let $W = \text{Span}\{\omega(a_p)-\omega(a_q) |p\neq q\}$. Then, for any $w\in W$ and $u\in w \mathbb{Q}\langle X \rangle$, $\widetilde{\textbf{I}_{g\circ \gamma}}[u]$ represents a convergent integral. The same argument as in (\ref{aux_4}) shows that $\textbf{A}[u] = \textbf{B}[u]$, so $u\in V$. \par
		In both cases, we constructed a codimension-one subspace $W$ of the words of weight $1$ such that $V$ is closed under left multiplication by $W$. If we can show that $V$ also contains all words of weight $1$, then Lemma \ref{co-dim1-shuffle} will complete the proof. This is done by choosing a suitable $B$: let $u=\Omega(g\gamma(1))$, and define $B$ by $$\widetilde{\textbf{I}_{g\circ \gamma}}[u] = B + \widetilde{\textbf{I}_{\gamma}}[u].$$
		
		Finally, for the case both $\gamma(0),\gamma(1)\in S$, split the path in half and apply the above case twice. 
	\end{proof}
	
	\section{Iterated integrals over a general base}
	\subsection{Definition and main properties}\label{MZVS_section}
	Let $\gamma$ be a path in the extended complex plane $\mathbb{P}^1 := \mathbb{C} \cup \{\infty\}$, and let $c_1,\cdots,c_n,d_1,\cdots,d_n \in \mathbb{P}^1$. Assume that $\gamma((0,1))$ (the image of the open interval under $\gamma$) does not contain any $c_i$ or $d_i$. Then
	\begin{equation}\label{iterint}\int_\gamma (\omega(c_1)-\omega(d_1))(\omega(c_2)-\omega(d_2))\cdots (\omega(c_n)-\omega(d_n))\qquad c_i,d_i\in \mathbb{P}^1\end{equation}
	converges if
	\begin{equation}\label{convcond}\gamma(0)\notin \{c_n,d_n\}, \qquad \gamma(1)\notin \{c_1,d_1\}.\end{equation}
	Indeed, if $\gamma$ lies entirely in $\mathbb{C}$ and $c_i, d_i\in \mathbb{C}$, this was noted above. If an endpoint is $\infty$, the integral still converges because $1/(x-a) - 1/(x-b) = O(1/x^2)$. \par
	
	Now we define the central object of our discussion:

	\begin{definition}
		Let $S$ be a finite subset of $\mathbb{P}^1$ and $n$ a positive integer. Define the $\mathbb{Q}$-vector space $W_n^{S}$ to be the $\mathbb{Q}$-span of all possible iterated integrals in (\ref{iterint}), where the $c_i,d_i$ range over all elements of $S$ and $\gamma$ ranges over all paths in $\mathbb{P}^1 - S$ with $$\gamma(0),\gamma(1)\in S,\qquad \gamma(0)\notin \{c_n,d_n\}, \qquad \gamma(1)\notin \{c_1,d_1\}.$$
		Set $$\MZV{S}{n} := \sum_k \sum_{i_1+\cdots+i_k = n} W_{i_1}^{S}\cdots W_{i_k}^{S}.$$
		We call $n$ the weight of $\MZV{S}{n}$. 
	\end{definition}
	
	We will soon see that $\MZV{S}{n}$ is a natural object of study. We shall always assume that $|S| \geq 3$. We first record some straightforward observations.
	\begin{proposition}\label{generalStrans}
		Assume that $|S| \geq 3$, and let $R$ be a rational function.
		\begin{enumerate}
			\item $2\pi i \in \MZV{S}{1}$
			\item $\MZV{S}{n} \MZV{S}{m} \subset \MZV{S}{n+m}$
			\item $\MZV{S}{n} \subset \MZV{R^{-1}(S)}{n}$
			\item If $R$ is invertible, then $\MZV{R(S)}{n} = \MZV{S}{n}$.
			\item If $R^{-1}(R(S))=S$, then $\MZV{R(S)}{n} \subset \MZV{S}{n}$.
		\end{enumerate}
	\end{proposition}
	\begin{proof}
		(1) Let $a,b,c$ be three distinct points of $S$, and consider the integral $\int_\gamma (\omega(b)-\omega(c))$, where $\gamma$ is a circle that starts and ends at $a$ and encloses $b$ but not $c$. The value of the integral is $\pm 2\pi i$, so this number is in $\MZV{S}{1}$.\par
		(2) This follows from the definition of $\MZV{S}{n}$ in terms of $W^{S}_i$.\par
		(3) It suffices to prove $W_n^{S} \subset W_n^{R^{-1}(S)}$. When $R$ is any rational function (not necessarily of degree $1$) we still have
		$$R^\ast \omega(a) = \frac{R'}{R-a}dx = \omega(R^{-1}(a)) - \omega(R^{-1}(\infty)) \qquad a\in \mathbb{P}^1,$$
		provided that we interpret the term $\omega(R^{-1}(a)):= \sum_{R(a_i) = a} \omega(a_i)$ counted with multiplicity.\par % This is because $R'/(R-a) = (R-a)'/(R-a)$, and for any rational function $f$, $f'/f$ has only simple poles, and residue of these poles corresponds to order of vanishing of $R$.
		Hence $$R^\ast (\omega(a) - \omega(b)) = \omega(R^{-1}(a)) - \omega(R^{-1}(b))$$ because the terms $\omega(R^{-1}(\infty))$ cancel. Let $a,b\in S$, and let $\gamma$ be any path from $a$ to $b$ for which \begin{equation}\label{convcond2}a\notin \{c_n,d_n\},\qquad b\notin \{ c_1,d_1\}.\end{equation}
		Consider the iterated integral $$I:=\int_\gamma (\omega(c_1)-\omega(d_1))(\omega(c_2)-\omega(d_2))\cdots (\omega(c_n)-\omega(d_n)) \in W^{S}_n.$$
		The map $R$ can be treated as a covering map away from its ramification points. Let $\widetilde{a}$ and $\widetilde{b}$ be any elements of $R^{-1}(a)$ and $R^{-1}(b)$, respectively. Then there exists a path $\widetilde{\gamma}$ (not necessarily unique) from $\widetilde{a}$ to $\widetilde{b}$ such that $R\circ \widetilde{\gamma} = \gamma$. Hence, by the pullback property of iterated integrals, 
		$$\begin{aligned}I &= \int_{R\circ \widetilde{\gamma}} (\omega(c_1)-\omega(d_1))(\omega(c_2)-\omega(d_2))\cdots (\omega(c_n)-\omega(d_n)) \\
			&= \int_{\widetilde{\gamma}} (\omega(R^{-1}(c_1))-\omega(R^{-1}(d_1)))\cdots (\omega(R^{-1}(c_n))-\omega(R^{-1}(d_n))).\end{aligned}$$
		This proves that $I\in W^{R^{-1}(S)}_{n}$. \par
		(4) follows from (3) by applying it to both $R$ and $R^{-1}$. (5) follows from (3) by replacing $S$ with $R(S)$.
	\end{proof}
	
	We record here a spanning set for $\MZV{S}{1}$. Recall that the cross-ratio of any four points $z_i \in \mathbb{P}^1$ is defined to be $\frac{(z_3-z_1)(z_4-z_2)}{(z_3-z_2)(z_4-z_1)}$; it is invariant under Möbius transformations.
	
	\begin{lemma}\label{level1span}
		The number $2\pi i$, together with the logarithms of the cross-ratios of all $4$-tuples of elements of $S$, spans $\MZV{S}{1}$.
	\end{lemma}
	\begin{proof}
		Follows from the formula
		$$\int_a^b (\omega(c)-\omega(d)) = \log \frac{(b-c)(a-d)}{(a-c)(b-d)}.$$
		The term $2\pi i$ arises from the branches of the logarithm.
	\end{proof}

	\begin{corollary}\label{corollaryofsupporttrans}
		Let $R$ be a rational function such that $R^{-1}(R(S)) = S$. Suppose that $\{0,1,\infty\}\subset R(S)$. Then
		$$\int_0^1 \omega(a_1)\cdots \omega(a_n) \in \MZV{S}{n}$$
		for $a_1\neq 1, a_n\neq 0, a_i\in R(S)$. Here the integration path can be any path from $0$ to $1$, not necessarily the straight line. 
	\end{corollary}
	\begin{proof}
		Because $R(S)$ contains $\{0,1,\infty\}$, we have
		$$\int_0^1 \omega(a_1)\cdots \omega(a_n) \in \MZV{R(S)}{n}.$$
		By Proposition \ref{generalStrans}, the right-hand side is contained in $\MZV{S}{n}$.
	\end{proof}
	
	Note that the condition $R^{-1}(R(S)) = S$ in the above corollary is automatic if $R$ is a Möbius transformation (i.e., invertible). The next theorem, which will be proved in the next section, is the central result of this paper. It states that $\CMZV{N}{n}$ is essentially $\MZV{S}{n}$ for certain $S$.
	
	\begin{theorem}[Main theorem]\label{CMZV_main_theorem}
		Let $N\geq 3$ and $S = \{0,\infty,1,\mu,\cdots,\mu^{N-1}\}$, where $\mu = e^{2\pi i / N}$. Then $\CMZV{N}{n} = \MZV{S}{n}$.
	\end{theorem}
	
	\subsection{Consequences of the Main Theorem \ref{CMZV_main_theorem}}
	
	The theorem has profound consequences for special values of multiple polylogarithms. We give a few examples below.
	
	\begin{corollary}\label{main_corollary}
		Let $N\geq 3$ and $S = \{0,\infty,1,\mu,\cdots,\mu^{N-1}\}$, where $\mu = e^{2\pi i / N}$. Let $R$ be a rational function such that $R^{-1}(R(S)) = S$ and $\{0,1,\infty\} \subset R(S)$. Then $$\int_0^1 \omega(a_1)\cdots \omega(a_n) \in \CMZV{N}{n}$$
		for $a_1\neq 1, a_n\neq 0, a_i\in R(S)$. 
	\end{corollary}
	\begin{proof}
		This follows from Corollary \ref{corollaryofsupporttrans} and the fact that $\MZV{S}{n} = \CMZV{N}{n}$.
	\end{proof}
	
	\begin{example}\label{level5Ex}
		Let $N=5$ and $\mu =e^{2\pi i /5}$. Let $R$ be the unique Möbius transformation such that $R^{-1}(0)=\mu$, $R^{-1}(1)=1$, and $R^{-1}(\infty)=\mu^2$. Then one checks that $R$ maps $S = \{0,\infty,1,\mu,\mu^2,\mu^3,\mu^4\}$ to
		$$R(S) = \left\{-\mu-\mu^2-\mu^3,1+\mu,1,0,\infty,\frac{\sqrt{5}+3}{2},\frac{\sqrt{5}+1}{2}\right\},$$
		Thus, when the $a_i$ are finite values in the above list, $\int_0^1 \omega(a_1)\cdots \omega(a_n)$ is a level-$5$ CMZV. Using (\ref{polylogtoint}), we see that, for the multiple polylogarithm $\Li_{s_1,\cdots,s_n}(x_1,\cdots,x_n)$, if $x_1^{-1},x_1^{-1}x_2^{-1},\cdots,x_1^{-1}\cdots x_n^{-1}$ are contained in $$\left\{-\mu-\mu^2-\mu^3,1+\mu,1,\frac{\sqrt{5}+3}{2} ,\frac{\sqrt{5}+1}{2} \right\},$$
		then $\Li_{s_1,\cdots,s_n}(x_1,\cdots,x_n) \in \CMZV{5}{n}$ with $n=\sum s_i$. In particular:
		\begin{itemize}
			\item $\Li_{s_1,\cdots,s_n}(z) \in \CMZV{5}{n}$ provided that $z = (\sqrt{5}-1)/2$ or $(3-\sqrt{5})/2$;
			\item $\Li_{s_1,\cdots,s_n}(x_1,\cdots,x_n) \in \CMZV{5}{n}$, provided that $n-2$ of the $x_i$ are equal to $1$ and the remaining two are equal to $(\sqrt{5}-1)/2$.
		\end{itemize}
	\end{example}
	
	\begin{example}\label{level6Ex}
		Let $N=6$ and $\mu =e^{2\pi i /6}$. Let $R$ be the Möbius transformation such that $R^{-1}(0)=\mu$, $R^{-1}(1)=1$, and $R^{-1}(\infty)=\mu^3$. Then $R$ maps $S = \{0,\infty,1,\mu,\mu^2,\mu^3,\mu^4,\mu^5\}$ to
		$$R(S)=\left\{1-i \sqrt{3},1+i \sqrt{3},1,0,-2,\infty,4,2\right\},$$
		Thus, when the $a_i$ are finite values in the above list, $\int_0^1 \omega(a_1)\cdots \omega(a_n)$ is a level-$6$ CMZV whenever it converges. For the multiple polylogarithm $\Li_{s_1,\cdots,s_n}(x_1,\cdots,x_n)$, if $x_1^{-1},x_1^{-1}x_2^{-1},\cdots,x_1^{-1}\cdots x_n^{-1}$ are contained in $$\left\{1-i \sqrt{3},1+i \sqrt{3},1,-2,4,2\right\},$$
		then $\Li_{s_1,\cdots,s_n}(x_1,\cdots,x_n) \in \CMZV{6}{n}$ with $n=\sum s_i$. In particular:
		\begin{itemize}
			\item The generalized polylogarithm $\Li_{s_1,\cdots,s_n}(z) \in \CMZV{6}{n}$ provided that $z = 1/2, -1/2$ or $1/4$;
			\item The multiple polylogarithm $\Li_{s_1,\cdots,s_n}(x_1,\cdots,x_n) \in \CMZV{6}{n}$ provided that $n-2$ of the $x_i$ are equal to $1$ and the remaining two are equal to $1/2$ or $-1/2$.
		\end{itemize}
	\end{example}
	
	\begin{example}\label{level6Ex_Deligne_value}
		Let $\mu = e^{2\pi i /6}$ and $S = \{0,1,\infty,\mu^2,\mu^4\}$. For $R(x) = \frac{x+1-\mu}{x-1}$, we have $R(S) = \{\mu^2,\infty,1,0,\mu\}$. Therefore,
		$$\int_0^1 \omega(a_1)\cdots \omega(a_n) \in \CMZV{3}, \qquad a_i\in \{0,1,\mu,\mu^2\}.$$
		Since $\mu$ is a primitive sixth root of unity, we see that certain values $\Li_{s_1,\cdots,s_n}(x_1,\cdots,x_n)$ for which the $x_i$ are sixth roots of unity actually lie in $\CMZV{3}{}$. 
	\end{example}
	
	\begin{example}\label{level10Ex1}
		Let $N=10$ and $\mu =e^{2\pi i /10}$. Let $R$ be the Möbius transformation such that $R^{-1}(0)=1$, $R^{-1}(1)=\mu^2$, and $R^{-1}(\infty)=\mu^6$. Then $R$ maps $S = \{0,\infty,1,\mu,\cdots,\mu^9\}$ to
		$$R(S) = \left\{\alpha,\bar{\alpha},0,\frac{1}{2},1,\frac{\sqrt{5}+1}{2},\frac{\sqrt{5}+3}{2} ,\sqrt{5}+3,\infty,-\sqrt{5}-2,\frac{-\sqrt{5}-1}{2},\frac{1-\sqrt{5}}{2}\right\},\qquad \alpha = \mu+\mu^3.$$
		Thus, when the $a_i$ are finite values in the above list, $\int_0^1 \omega(a_1)\cdots \omega(a_n)$ is a level-$10$ CMZV whenever it converges. 
		When $\frac{1}{2}$ is present among the $a_i$, the integration path for $\int_0^1 \omega(a_1)\cdots \omega(a_n)$ must be deformed to avoid this point. For any such deformation, the assertion that the integral belongs to $\CMZV{10}{}$ remains valid. 
	\end{example}

	In the examples above, we used only the case in which $R$ is invertible. We next consider some higher-degree rational functions $R$.\par
	
	%Let $G$ be the group of Möbius that sends a finite set $S$ to $S$. For each subgroup $H$ of $G$, and $R_0$ Möbius, $R = \sum_{\sigma\in H} R_0\circ \sigma$ is $H$-invariant, and $\deg R \leq |H|$. Hence if $R$ is not a constant, for $s\in \mathbb{P}^1$, the fibre of $R^{-1}(R(s))$ are exactly $Hs$, counted with multiplicity. In particular, $R$ satisfies $R^{-1}(R(S))\subset S$. \par
	
	\begin{example}
		Let $N=10$, $\mu =e^{2\pi i /10}$, $S=\{0,\infty,1,\cdots,\mu^9\}$, and $R(x) = x+\mu/x$. The set $R(S)$ has $6$ elements: $\{\infty, R(1), R(\mu^2),R(\mu^3),R(\mu^4),R(\mu^5)\}$, and $R^{-1}(R(S)) = S$. For each $3$-tuple of this set, choose a Möbius transformation $R_1$ that maps the tuple to $(0,1,\infty)$. As an illustration, consider the $R_1$ such that $R_1(R(\mu),R(\mu^2),0) = (0,1,\infty)$. Then 
		$$R_1(R(S)) = \left\{\frac{-\sqrt{5}-1}{2},0,1,\infty,-\sqrt{5}-2,-\sqrt{5}-1\right\}.$$
		Thus, when the $a_i$ are finite values in the above list, $\int_0^1 \omega(a_1)\cdots \omega(a_n)$ is a level-$10$ CMZV. This cannot be proved using degree-one rational functions $R$ alone (as in the previous examples). \par
		If we choose another $3$-tuple that maps to $(0,1,\infty)$, say $R_2(R(\mu^4),R(\mu^2),R(\mu^3)) = (0,1,\infty)$, then 
		$$R_2(R(S)) = \left\{\frac{1}{2},\frac{\sqrt{5}+1}{4} ,1,\infty,0,\frac{3-\sqrt{5}}{4} \right\}.$$
		Thus, when the $a_i$ are finite values in the above list, $\int_0^1 \omega(a_1)\cdots \omega(a_n)$ is a level-$10$ CMZV. 
	\end{example}
	
	\begin{proposition}\label{zinpolylog}
		Let $\mu = e^{2\pi i /N}$, where $N\geq 3$, and let $S=\{0,\infty,1,\mu, \cdots,\mu^{N-1}\}$. Let $R$ be a rational function such that $R^{-1}(R(S))= S$ and $\{0,1,\infty\}\subset R(S)$. Then, for any $z\in R(S) - \{0,1,\infty\}$, $$\Li_{s_1,\cdots,s_n}(z) \in \CMZV{N}{s_1+\cdots+s_n}.$$
	\end{proposition}
	\begin{proof}
		Writing $n=s_1+\cdots+s_n$, we have $$\Li_{s_1,\cdots,s_n}(z) = (-1)^n \int_0^z \omega(0)^{s_1-1}\omega(1) \omega(0)^{s_2-1}\omega(1)\cdots \omega(0)^{s_n-1} \omega(1) \in \MZV{\{0,1,\infty,z\}}{n}\subset \MZV{R(S)}{n}.$$
		Since $R^{-1}(R(S)) = S$, Corollary \ref{main_corollary} implies that the above space is contained in $\MZV{S}{n}$, which equals $\CMZV{N}{n}$ by the main theorem.
	\end{proof}
	
	\begin{example}\label{level10Ex2}
		Let $N=10$, $\mu =e^{2\pi i /10}$, and $S=\{0,\infty,1,\cdots,\mu^9\}$. Let $R(x) = x+1/x$. One checks that $R^{-1}(R(S))\subset S$. The set $R(S)$ has $7$ elements. For each $3$-tuple of $R(S)$, choose a Möbius transformation $R_1$ that maps the tuple to $(0,1,\infty)$. For example, if
		$R_1(R(1),R(\mu^3),R(\mu)) = (0,1,\infty)$, then 
		$$R_1(R(S)) = \left\{\frac{3 \sqrt{5}-5}{2},0,\infty,5 \sqrt{5}-10,1,\frac{15-5 \sqrt{5}}{4} ,4 \sqrt{5}-8\right\}.$$ 
		If $R_2(R(\mu^5),R(\mu^2),R(1)) = (0,1,\infty)$, then 
		$$R_2(R(S)) = \left\{2 \sqrt{5}-5,\infty,5,1,45-20 \sqrt{5},9-4 \sqrt{5},0\right\}.$$
		Consequently, the generalized polylogarithm $\Li_{s_1,\cdots,s_n}(z)$ is an element of $\CMZV{10}{s_1+\cdots+s_n}$ whenever $$z=2 \sqrt{5}-5,\quad \frac{1}{5},\quad 45-20 \sqrt{5},\quad 9-4 \sqrt{5},\quad \frac{15-5 \sqrt{5}}{4},\quad \text{or}\quad4 \sqrt{5}-8.$$
	\end{example}
	
	We can generate many more examples. Nonetheless, for a fixed level $N$, the number of possibilities is finite. 
	
	\begin{proposition}\label{S_unit_eq_finite} Let $S$ be any finite subset of $\mathbb{P}^1$ with $|S|\geq 3$. \\
		(a) The number of rational functions $R$ such that $$R^{-1}(0), R^{-1}(1), R^{-1}(\infty) \subset S$$ is finite. \\
		(b) The number of rational functions $R$ such that $\{0,1,\infty\} \subset R(S)$ and $R^{-1}(R(S)) = S$ is finite. 
	\end{proposition}
	\begin{proof}
		(a) The displayed condition is equivalent to the fact that the divisors of $R$ and $1-R$ are supported on $S$. This is an \textit{$S$-unit equation} over a genus-$0$ function field \cite{silverman1984s, mason1983hyperelliptic, brunault2020k_4}, and it is known to have only finitely many solutions, all of which satisfy $\deg(R)\leq |S|-2$. \\
		(b) Any such $R$ satisfies the condition in (a).
	\end{proof}

	\section{Proof of the Main Theorem \ref{CMZV_main_theorem}}
	We first develop two key theorems that are valid for a general finite set $S\subset \mathbb{P}^1$. Since $\MZV{S}{n} = \MZV{R(S)}{n}$ for any Möbius transformation $R$, we may assume that $\infty \in S$. The space $W^S_n$ is the $\mathbb{Q}$-span of iterated integrals of the form 
	\begin{equation}\label{aux_7}\int_\gamma \omega(c_1)\cdots \omega(c_n), \qquad c_1\neq \gamma(1), c_n\neq \gamma(0).\end{equation}
	
	For the proofs of the next two theorems, we assume that $S = \{\infty,a_1,\cdots,a_k\}$ and let $X = \{\omega(a_1),\cdots,\omega(a_k)\}$. We use the notation introduced at the beginning of Section \ref{tangential_base_point_section}.
	
	\begin{theorem}\label{homotopymodzero}
		In the setting of (\ref{aux_7}), if $\gamma$ is a loop, then, for $n\geq 2$, the iterated integral is in $$\overline{\MZV{S}{n}} := \sum_{1\leq k < n} \MZV{S}{k}\MZV{S}{n-k},$$
		i.e., it is a linear combination of products of elements of strictly lower weights.
	\end{theorem}
	\begin{proof}
		We can perturb $\gamma$ by $\varepsilon>0$ to obtain a path $\gamma(\varepsilon)$ based at $a_1+\varepsilon$. Consider the homomorphism $$\pi_1(\mathbb{P}^1-S, a_1+\varepsilon) \to \mathbb{C}\llangle X \rrangle, \qquad p \mapsto \sum_{\omega\in X^\ast}\left(\int_{p} \omega \right)\omega = \textbf{I}_{p}.$$
		The group $\pi_1(\mathbb{P}^1-S, a_1+\varepsilon)$ is free of rank $|S|-1 = k$. For each $1\leq i\leq k$, a generator can be viewed as a loop based at $a_1+\varepsilon$ that encloses only $a_i$; we denote this generator by $\gamma_i(\varepsilon)$. We first find an explicit expression for each $\textbf{I}_{\gamma_i(\varepsilon)}$.
		
		For $i=1$, Lemma \ref{radius0} implies that $\textbf{I}_{\gamma_1(\varepsilon)} = \exp(2\pi i \omega(a_1)) + O(\varepsilon)$. We next investigate the other values of $i$. Without loss of generality, we focus on $i=2$ and deform $\gamma_2(\varepsilon)$ into the following path:
		\begin{figure}[h]
			\centering
			\begin{tikzpicture}[decoration={markings,
					%mark=at position 0.6cm with {\arrow[line width=1pt]{>}},
					mark=at position 5cm with {\arrow[line width=1pt]{>}},
					mark=at position 8.85cm with {\arrow[line width=1pt]{>}},
					mark=at position 13cm with {\arrow[line width=1pt]{>}}
				}
				]
				
				% The path
				\path[draw,line width=0.8pt,postaction=decorate] (-0.3,-0.15) -- (-0.3,0.15) -- (7.43,0.15)  \centerarcpath(8,0)(165.522:-165.522:0.6) -- (-0.3,-0.15);
				
				\filldraw[black] (-0.3,0) circle (1pt) node[anchor=west]{$a_1+\varepsilon$};
				\filldraw[black] (-0.8,0) circle (1pt) node[anchor=east]{$a_1$};
				\filldraw[black] (8,0) circle (1pt) node[anchor=east]{$a_2$};
				\node at (4,0.4) {$\Gamma(\varepsilon)$};
				\node at (8.7,0.4) {$C$};
				\end{tikzpicture}\caption{Integration along a loop enclosing $a_2$ and based at $a_1+\varepsilon$.}
		\end{figure}
		
		As in the figure, $$\textbf{I}_{\gamma_2(\varepsilon)}:= \sum_{\omega\in X^\ast} \int_{\gamma_2(\varepsilon)} \omega = \textbf{I}_{\Gamma(\varepsilon)}^{-1} \textbf{I}_{C} \textbf{I}_{\Gamma(\varepsilon)}.$$
		From Proposition \ref{mainregtheorem}, there exists $A_{2,\varepsilon}, B_{2,\varepsilon} \in \mathbb{C}$ such that
		$$\textbf{I}_{\Gamma(\varepsilon)} = e^{A_{2,\varepsilon} \omega(a_2)} \widetilde{\textbf{I}_2} e^{B_{2,\varepsilon} \omega(a_1)} + O(\varepsilon^{1/2}),$$
		where $\widetilde{\textbf{I}_2}$ is the unique group-like lift of $\lim_{\varepsilon\to 0} \pi(\textbf{I}_{\Gamma(\varepsilon)})$ for which the coefficients of $\omega(a_1)$ and $\omega(a_2)$ are $0$, and $$B_{2,\varepsilon} = \int_{\Gamma(\varepsilon)} \omega(a_1) = \int_{a_1+\varepsilon}^{a_2-\varepsilon} \frac{1}{x-a_1} dx = \log(a_2-a_1)-\log\varepsilon + O(\varepsilon^{1/2}).$$ The coefficients of $\widetilde{\textbf{I}_2}$ lie in $\MZV{S}{n}$: for coefficients arising from convergent integrals, this follows from the definition of $\MZV{S}{n}$; for divergent words, it follows from the recurrence (\ref{reg}). Moreover, Lemma \ref{radius0} implies that $\textbf{I}_{C} = e^{-2\pi i \omega(a_2)} + O(\varepsilon)$, and hence
		$$\begin{aligned}\textbf{I}_{\gamma_2(\varepsilon)} &= \left( e^{A_{2,\varepsilon} \omega(a_2)} \widetilde{\textbf{I}_2} e^{B_{2,\varepsilon} \omega(a_1)} \right)^{-1} e^{- 2\pi i \omega(a_2)}   \left( e^{A_{2,\varepsilon} \omega(a_2)} \widetilde{\textbf{I}_2} e^{B_{2,\varepsilon} \omega(a_1)} \right) + O(\varepsilon^{1/2})  \\ 
			&= e^{-B_{2,\varepsilon} \omega(a_1)} \widetilde{\textbf{I}_2}^{-1}  e^{- 2\pi i \omega(a_2)}\widetilde{\textbf{I}_2} e^{B_{2,\varepsilon} \omega(a_1)} + O(\varepsilon^{1/2})
		\end{aligned}.$$
		Therefore, for each generator $\gamma_j(\varepsilon)$ of the fundamental group $\pi_1(\mathbb{P}^1-S,a_1+\varepsilon)$, we have proved that
		$$\textbf{I}_{\gamma_j(\varepsilon)} = \begin{cases}e^{2\pi i \omega(a_1)} + O(\varepsilon)\quad &j=1, \\ e^{-B_{j,\varepsilon} \omega(a_1)} \widetilde{\textbf{I}_j}^{-1}  e^{- 2\pi i \omega(a_j)}\widetilde{\textbf{I}_j} e^{B_{j,\varepsilon} \omega(a_1)} + O(\varepsilon^{1/2}) \quad &j\geq 2,\\ \end{cases}$$
		where $B_{j,\varepsilon} = \log(a_j-a_1)-\log\varepsilon + O(\varepsilon^{1/2})$ and the weight-$n$ coefficients of $\widetilde{\textbf{I}_j}$ are in $\MZV{S}{n}$.
		
		We need to show that, for any word $u$ of weight $n$ that neither starts nor ends with $\omega(a_1)$ and for any product $\textbf{P}$ of $\{\textbf{I}_{\gamma_1(\varepsilon)},\textbf{I}_{\gamma_2(\varepsilon)},\cdots,\textbf{I}_{\gamma_k(\varepsilon)}\}$, we have \begin{equation}\label{aux_6}\lim_{\varepsilon\to 0}\textbf{P}[u] \in \overline{\MZV{S}{n}}.\end{equation}
		
		Observe that, for arbitrary $\textbf{A}, \textbf{B}\in \mathbb{C}\llangle X \rrangle$ whose coefficients of weight-$n$ monomials lie in $\MZV{S}{n}$ and for any $u\in X^\ast$, one has $$\textbf{A}\textbf{B}[u] \equiv \textbf{A}[u] + \textbf{B}[u] \pmod{\overline{\MZV{S}{n}}}, \qquad \textbf{A}^{-1}[u] = -\textbf{A}[u] \pmod{\overline{\MZV{S}{n}}}.$$
		
		By this observation, it suffices to show that $B_{j,\varepsilon} - B_{i,\varepsilon} \in \MZV{S}{1} + O(\varepsilon)$ for $i,j\neq 1$: the term $e^{(B_{j,\varepsilon} - B_{i,\varepsilon})\omega(a_1)}$ occurs in the product $\textbf{I}_{\gamma_i(\varepsilon)} \textbf{I}_{\gamma_j(\varepsilon)}$. This holds because
		$$B_{j,\varepsilon} - B_{i,\varepsilon} = \log\frac{a_j-a_1}{a_i-a_1}+ O(\varepsilon),$$
		This belongs to $\MZV{S}{1}$ by Lemma \ref{level1span}, because $\frac{a_j-a_1}{a_i-a_1}$ is the cross-ratio of the four points $\{\infty,a_j,a_i,a_1\} \subset S$. 
	\end{proof}
	
	For our given finite set $S$, let $G$ denote its symmetry group: $$G:= \{ \text{Möbius transformation } R \: | \: R(S) = S\}.$$ We will see that $G$ has a significant effect on the structure of $\MZV{S}{n}$. \par
	\begin{definition}
		Let $T = \{(t_1,s_1),\cdots,(t_r,s_r)\}$ be a subset of $S^2$. We define an undirected graph $\mathcal{G}(S,T)$ as follows: the vertex set is $S$, and for each $(t_i,s_i)\in T$ and each $g\in G$, connect $gt_i$ and $gs_i$ with an edge. We call $T$ a set of \textit{complete edges} if the graph $\mathcal{G}(S,T)$ is connected.
	\end{definition}
	
	\begin{example}\label{complete_edge_ex}
		(a) Let $\mu = e^{2\pi i /N}$ and $S = \{0,\infty,1,\mu,\cdots,\mu^{N-1}\}$. Then $G$ contains $x\mapsto 1/x$ and $x\mapsto \mu x$. We claim that $T = \{(0,1)\}$ is a set of complete edges. Indeed, by repeatedly applying $x\mapsto \mu x$, the vertices $1,\mu,\cdots,\mu^{N-1}$ lie in the connected component of $0$. Applying $x\mapsto 1/x$ to $(0,1)$ gives $(\infty, 1)$, so $\infty$ also lies in the connected component of $0$. Thus, $\mathcal{G}(S,T)$ is connected.\\
		(b) Viewing $\mathbb{P}^1$ as the unit sphere, let $S \subset \mathbb{P}^1$ be the set of vertices of a Platonic solid (or, more generally, an Archimedean solid). Then any single edge is a set of complete edges. 
	\end{example}
	
	We return to the setting of a general $S$. Let $T = \{(t_1,s_1),\cdots,(t_r,s_r)\}$ be a subset of $S^2$, and choose fixed paths $\gamma_1,\cdots,\gamma_r$ whose endpoints are $t_i$ and $s_i$, respectively. Let $W^{S,T}_{n}$ be\footnote{It would be more precise to include $\gamma_1,\cdots,\gamma_r$ in the notation, but we do not do so.} the $\mathbb{Q}$-space spanned by 
	$$\int_{\gamma_i} \omega(c_1)\cdots \omega(c_n), \qquad c_1\neq \gamma_i(1), c_n\neq \gamma_i(0), \qquad 1\leq i\leq r,c_i\in \{a_1,\cdots,a_k\} = S - \{\infty\}.$$
	Define $$\MZV{S,T}{n} := \sum_k \sum_{i_1+\cdots+i_k = n} W^{S,T}_{i_1}\cdots W^{S,T}_{i_k}.$$
	\begin{theorem}\label{graphconnectedtheo}
		If $T$ is a set of complete edges, then\footnote{Here, the tensor product is graded by weight, with $\MZV{S}{1}$ assigned weight $1$.} $\MZV{S,T}{n} \otimes_\mathbb{Q} \MZV{S}{1} = \MZV{S}{n}$.
	\end{theorem}
	In essence, this result says that, in the definition of $\MZV{S}{n}$, instead of allowing $\gamma$ to range over all possible paths (of which there are infinitely many), it suffices to take $\gamma\in \{\gamma_1,\cdots,\gamma_r\}$.
	\begin{proof}
		Let $\widetilde{\textbf{I}_{\gamma_i}}$ be the regularized series given by Proposition \ref{mainregtheorem}.\footnote{We can always parametrize the path $\gamma_i$ so that it is regular at both endpoints.} Its weight-$n$ coefficients lie in $\MZV{S,T}{n}$: for convergent integrals, this is simply the definition of $\MZV{S,T}{n}$; for divergent integrals, it follows from the recurrence (\ref{reg}). We prove the statement by induction on $n$. The case $n=1$ holds because we have already tensored with the weight-$1$ space. Now assume that $n\geq 2$. \par 
		Let $\gamma$ be an arbitrary path as in (\ref{aux_7}). By the definition of completeness for $T$, there exists a path in the graph $\mathcal{G}(S,T)$ connecting $\gamma(0)$ and $\gamma(1)$, say $$(g_1t_1,g_1s_1), (g_2t_2,g_2s_2), \cdots, (g_qt_q,g_qs_q),$$
		with $g_i s_i = g_{i+1} t_{i+1} := b_i$, $g_1t_1 = \gamma(0) := b_0$, and $g_qs_q = \gamma(1) = b_q$. By removing a cycle if necessary, we may assume that the above path in $\mathcal{G}(S,T)$ contains no cycles, i.e., that $b_0,\cdots,b_q$ are pairwise distinct. \par
		Write $\iota := (g_q\gamma_q) \cdots (g_1\gamma_1)$. We slightly deform $\iota$ into a path $\iota(\varepsilon)$ that does not pass through any point of $S$, as shown in the figure. 
		\begin{figure}[h]
			\centering
			\includegraphics[scale=0.7]{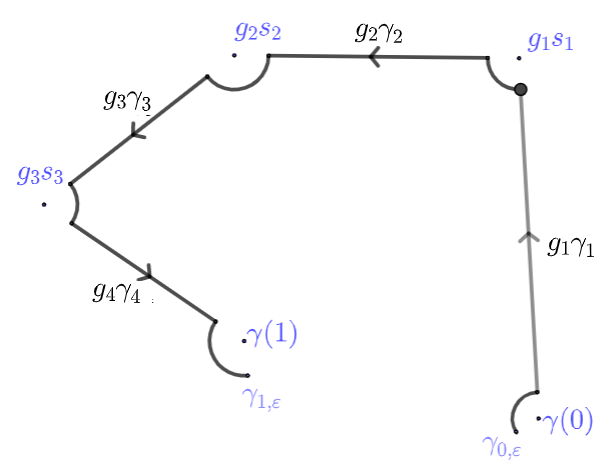}
			\caption{\small The path $\iota(\varepsilon)$. Here the paths $g_i\gamma_i$ are restricted to the interval $[\varepsilon,1-\varepsilon]$.}
		\end{figure}
		Also, let $\rho := \gamma^{-1}\iota$; it is a loop based at $\gamma(0)$. We perturb $\rho$ slightly into $\rho(\varepsilon)$ so that it does not pass through any point of $S$ and is a loop based at $\gamma_{0,\varepsilon}$. Let $\gamma(\varepsilon) := \iota(\varepsilon) \rho(\varepsilon)$. Our goal is to show that, for $u\in X^\ast$ that does not start with $\omega(\gamma(1))$ or end with $\omega(\gamma(0))$, we have
		$$\lim_{\varepsilon\to 0}\textbf{I}_{\gamma(\varepsilon)}[u] \in \MZV{S,T}{n}\otimes \MZV{S}{1}.$$
		We begin by finding more explicit expressions for $\textbf{I}_{\iota(\varepsilon)}$ and $\textbf{I}_{\rho(\varepsilon)}$. Let $$\textbf{J}_i := \sum_{w\in X^\ast} \left(\int_{g_i\gamma_i | [\varepsilon,1-\varepsilon]} w \right)w \in \mathbb{C}\llangle X \rrangle$$ be the element associated with $g_i \gamma_i$ as shown in the figure, and let $\textbf{C}_i\in \mathbb{C}\llangle X \rrangle$ be the element associated with the circular arcs. We have $\textbf{I}_{\iota(\varepsilon)} = \textbf{C}_q \textbf{J}_q \cdots \textbf{C}_1 \textbf{J}_1 \textbf{C}_0$. For the element $\widetilde{\textbf{J}_i}$ defined by Proposition \ref{mainregtheorem}, we have
		$$\textbf{J}_i = e^{(B_i + \log \varepsilon)\Omega(b_i)} \widetilde{\textbf{J}_i }e^{(A_i - \log \varepsilon)\Omega(b_{i-1})} + O(\varepsilon^{1/2})$$
		for some $A_i, B_i\in \mathbb{C}$. Moreover, Lemmas \ref{radius0} and \ref{radiusinf} imply that
		$$\textbf{C}_i = e^{C_i \Omega(b_i)} + O(\varepsilon),$$
		for some $C_i \in \mathbb{C}$. Substituting these expressions into $\textbf{I}_{\iota(\varepsilon)} = \textbf{C}_q \textbf{J}_q \cdots \textbf{C}_1 \textbf{J}_1 \textbf{C}_0$ and renaming the constants yields
		\begin{equation}\label{aux_8}\textbf{I}_{
				\iota(\varepsilon)} = e^{B_\varepsilon \Omega(b_q)} \widetilde{\textbf{J}}_q e^{A_q \Omega(b_q)} \widetilde{\textbf{J}}_{q-1} e^{A_{q-1} \Omega(b_{q-1})}\cdots \widetilde{\textbf{J}}_{2} e^{A_2 \Omega(b_1)} \widetilde{\textbf{J}}_{1} e^{C_\varepsilon \Omega(b_0)} + O(\varepsilon^{1/2})\end{equation}
		for some $A_i\in \mathbb{C}$ independent of $\varepsilon$. The element $\widetilde{\textbf{J}}_{i}$ is a regularized integral along the path $g_i\gamma_i$, so Lemma \ref{group_trans_reg} implies $\widetilde{\textbf{J}}_{i} = e^{D_i \Omega(b_i)} (g_i \cdot \widetilde{\textbf{I}_i}) e^{E_i \Omega(b_{i-1})}$
		for some constants $D_i, E_i$. Substituting this expression above and renaming the constants, we obtain
		$$\textbf{I}_{
			\iota(\varepsilon)} = e^{B_\varepsilon \Omega(b_q)} (g_q\cdot \widetilde{\textbf{I}}_{q}) e^{A_q \Omega(b_q)}(g_{q-1}\cdot \widetilde{\textbf{I}}_{q-1}) e^{A_{q-1} \Omega(b_{q-1})}\cdots (g_2\cdot \widetilde{\textbf{I}}_{2}) e^{A_2 \Omega(b_1)} (g_1\cdot \widetilde{\textbf{I}}_{1}) e^{C_\varepsilon \Omega(b_0)} + O(\varepsilon^{1/2}),$$
		where $B_\varepsilon = B+\log\varepsilon$ and $C_\varepsilon = C-\log\varepsilon$ for some $B,C\in \mathbb{C}$. We now turn to $\textbf{I}_{\rho(\varepsilon)}$. By the previous lemma,
		$$\textbf{I}_{\rho(\epsilon)} = e^{F_\varepsilon \Omega(b_0)} \textbf{M} e^{G_\varepsilon \Omega(b_0)},$$
		where $F_\varepsilon = F + \log \varepsilon$ and $G_\varepsilon = G-\log \varepsilon$ for some $F,G\in \mathbb{C}$, and $\textbf{M}$ is a formal power series whose weight-$n$ coefficients are in $\overline{\MZV{S}{n}}$.
		Therefore, again renaming constants, \begin{equation}\label{aux_9}\textbf{I}_{\gamma(\epsilon)} = \textbf{I}_{\iota(\epsilon)}\textbf{I}_{\rho(\epsilon)} = e^{B_\varepsilon \Omega(b_q)} (g_q\cdot \widetilde{\textbf{I}}_{q}) e^{A_q \Omega(b_q)}\cdot \cdots (g_2\cdot \widetilde{\textbf{I}}_{2}) e^{A_2 \Omega(b_1)} (g_1\cdot \widetilde{\textbf{I}}_{1}) e^{A_1 \Omega(b_0)} \textbf{M} e^{G_\varepsilon \Omega(b_0)} + O(\varepsilon^{1/2}).\end{equation}
		Recall that our goal is to show that $\lim_{\varepsilon\to 0} \textbf{I}_{\gamma(\varepsilon)}[u] \in \MZV{S,T}{n}\otimes \MZV{S}{1}$. The exponentials at the front and end can be ignored. By the induction hypothesis on $n$, we may assume that the weight-$n$ coefficients of $\textbf{M}$ lie in this space. We saw at the beginning of the proof that the coefficients of $\widetilde{\textbf{I}}_{i}$ lie in $\MZV{S,T}{n}$, and therefore so do those of $g_i\cdot \widetilde{\textbf{I}}_{i}$. It remains to show that $A_1,\cdots,A_q\in \MZV{S}{1}$. Recall that $b_0,b_1,\cdots,b_q$ are pairwise distinct elements of $S$. \par
		First, assume that each $b_i$ is finite, so $\Omega(b_i) = \omega(b_i)$. For each $i$, comparing the coefficients of $\omega(b_i)$ on both sides of (\ref{aux_9}) gives
		$$\int_\gamma \omega(b_i) = A_i + \sum_j \widetilde{\textbf{I}}_j[\omega(b_i)].$$
		The left-hand side is in $\MZV{S}{1}$. If a term $\widetilde{\textbf{I}}_j[\omega(b_i)]$ on the right-hand side arises from a convergent integral, then it is in $\MZV{S}{1}$; otherwise, by our choice of $\widetilde{\textbf{I}}_j$, it is $0$. Therefore, $A_i \in \MZV{S}{1}$. This completes the proof when each $g_it_i$ is finite. \par 
		Finally, when one of the $b_i$ equals $\infty$, we have $\Omega(b_i) = -\omega(a_1)-\cdots-\omega(a_k)$. The argument is largely parallel to the one above: we can still solve for $A_1,\cdots,A_q$ and conclude that they are in $\MZV{S}{1}$.
	\end{proof}
	
	\begin{lemma}
		Let $S = \{0,\infty,1,e^{2\pi i /N},\cdots,e^{2\pi i (N-1)/N}\}$. When $N\geq 3$, $\CMZV{N}{1} = \MZV{S}{1}$. 
	\end{lemma}
	\begin{proof}
		The containment $\CMZV{N}{1}\subset \MZV{S}{1}$ is evident. For the reverse inclusion, Lemma \ref{level1span} shows that it suffices to prove that the logarithm of the cross-ratio of each $4$-tuple in $S$ can be written as a linear combination of $2\pi i$ and $\log(1-\mu^i)$, where $\mu=e^{2\pi i /N}$. We omit this straightforward computation.
	\end{proof}
	
	Finally, we can prove our main result.
	
	\begin{proof}[Proof of Theorem \ref{CMZV_main_theorem}]
		Let $S = \{0,\infty,1,e^{2\pi i /N},\cdots,e^{2\pi i (N-1)/N}\}$. From Example \ref{complete_edge_ex}, we know that $T=\{(0,1)\}$ is a set of complete edges for $S$. Choosing $\gamma$ to be the straight line from $0$ to $1$, the corresponding space $\MZV{S,T}{n}$ is the $\mathbb{Q}$-span of the numbers
		$$\int_\gamma \omega(c_1)\cdots \omega(c_n) \qquad c_i\in \{0,1,e^{2\pi i /N},\cdots,e^{2\pi i (N-1)/N}\}, c_1\neq 1, c_n\neq 0,$$
		so $\MZV{S,T}{n}=\CMZV{N}{n}$. Theorem \ref{graphconnectedtheo} then implies that $\CMZV{N}{n}\otimes_\mathbb{Q} \MZV{S}{1} = \MZV{S}{n}$. The previous lemma shows that the weight-$1$ component can be absorbed.
	\end{proof}
	
	\begin{remark}
		The equality $\CMZV{N}{n}\otimes_\mathbb{Q} \MZV{S}{1} = \MZV{S}{n}$ holds for every level $N$. However, when $N=1$ or $2$, we cannot absorb the weight-$1$ space; in these two cases, we have $$\CMZV{N}{n}\otimes_\mathbb{Q} \mathbb{Q}[2\pi i] = \MZV{S}{n}, \qquad N\in \{1,2\}.$$
	\end{remark}
	
	\subsection{Explicit computations}\label{expcomp}
	Our proof of Theorem \ref{CMZV_main_theorem} provides a way to make the inclusion in Corollary \ref{main_corollary} explicit. We first illustrate this with an example. For a fixed level $N$, choose a standard group-like element $\widetilde{\textbf{I}} \in \mathbb{C}\llangle X \rrangle$, where $X = \{\omega(0),\omega(1),\omega(e^{2\pi i/N}),\cdots,\omega(e^{2\pi i (N-1)/N})\}$. For $w\in X^\ast$ that does not start with $\omega(1)$ or end with $\omega(0)$, define $\widetilde{\textbf{I}}[w] = \int_0^1 w$. Then $\widetilde{\textbf{I}}$ is the unique group-like element satisfying $\widetilde{\textbf{I}}[\omega(1)] = \widetilde{\textbf{I}}[\omega(0)]=0$, and its weight-$n$ coefficients are in $\CMZV{N}{n}$. 
	
	\begin{example}Let us revisit Example \ref{level5Ex}. Write $S=\{0,\infty,1,\cdots,\mu^4\}$, $\mu = e^{2\pi i /5}$, and $X=\{\omega(0),\omega(1),\omega(\mu),\cdots,\omega(\mu^4)\}$. There, we used the Möbius transformation $R$ that maps $(\mu,1,\mu^2)$ to $(0,1,\infty)$ to assert that
		$$\int_0^1 \omega(c_1)\cdots\omega(c_n) \in \CMZV{5}{n} \quad \text{whenever }c_i \in \left\{-\mu-\mu^2-\mu^3,1+\mu,1,0,\frac{\sqrt{5}+3}{2},\frac{\sqrt{5}+1}{2}\right\} = R(S) - \{\infty\}.$$
		
		How can the integral be expressed explicitly as an element of $\CMZV{5}{n}$? Let $[0,1]$ denote the straight-line path from $0$ to $1$. Explicitly, we have $R(x) = (1 + \mu) (\mu - x)/(\mu^2 - x)$ and $$\int_0^1 \omega(c_1)\cdots\omega(c_n) = \int_{R^{-1}[0,1]} R^\ast\omega(c_1)\cdots R^\ast\omega(c_n).$$ Note that $R^{-1}[0,1]$ is a circular arc from $\mu$ to $1$. 
		
		\begin{figure}[h]
			\centering
			\begin{tikzpicture}[scale=5, line width=1pt]
				\draw[help lines,->] (-0.2,0) -- (1.2,0) coordinate (xaxis);
				\draw[help lines,->] (0,-0.2) -- (0,1) coordinate (yaxis);
				\draw [shift={(0,0)}]  plot[domain=0:1.1949282492261746,variable=\t]({1*0.09436863320044009*cos(\t r)+0*0.09436863320044009*sin(\t r)},{0*0.09436863320044009*cos(\t r)+1*0.09436863320044009*sin(\t r)});
				\draw [shift={(1,0)}]  plot[domain=1.6018297099869439:3.141592653589793,variable=\t]({1*0.06205680441467056*cos(\t r)+0*0.06205680441467056*sin(\t r)},{0*0.06205680441467056*cos(\t r)+1*0.06205680441467056*sin(\t r)});
				\draw [shift={(0.3670801437747959,0.9301893183896895)}]  plot[domain=4.336520902815968:5.8703352300664084,variable=\t]({1*0.0739471405391017*cos(\t r)+0*0.0739471405391017*sin(\t r)},{0*0.0739471405391017*cos(\t r)+1*0.0739471405391017*sin(\t r)});
				\draw [middlearrow={>}] (0.33993561679396744,0.861404478134756)-- (0.03464085144304853,0.087780694594084);
				\draw [middlearrow={>}] (0.09436863320044009,0)-- (0.937943195585329,0);
				\draw [blue,middlearrow={<}, shift={(0,0)}]  plot[domain=0.06206676638409518:1.1209642501372628,variable=\t]({1*1*cos(\t r)+0*1*sin(\t r)},{0*1*cos(\t r)+1*1*sin(\t r)});
				\filldraw[black] (0,0) circle (0.2pt) node[anchor=east]{$0$};
				\filldraw[black] (1,0) circle (0.2pt) node[anchor=west]{$1$};
				\filldraw[black] (0.37,0.93) circle (0.2pt) node[anchor=south]{$e^{2\pi i /5}$};
				\node[blue] at (1.05,0.5) {$R^{-1} [0,1]$};
				\node at (0.5,-0.1) {$[0,1]$};
				\node at (0.1,0.6) {$g[0,1]^{-1}$};
			\end{tikzpicture}\caption{\small Deforming the integration path $R^{-1}[0,1]$ to two parts: $g[0,1]^{-1}$ and $[0,1]$.}
		\end{figure}
		
		The integration path $R^{-1}[0,1]$ can be deformed into two segments: the path $g[0,1]^{-1}$ followed by $[0,1]$. Here $g\in G$ is rotation through an angle of $2\pi/5$. For $\textbf{I}_{R^{-1}[\varepsilon,1-\varepsilon]}$ as given in Proposition \ref{mainregtheorem}, we have
		$$\textbf{I}_{R^{-1}[\varepsilon,1-\varepsilon]} = e^{A_\varepsilon \omega(1)} \widetilde{\textbf{I}} e^{B\omega(0)} (g\cdot \widetilde{\textbf{I}}^{-1}) e^{C_\varepsilon \omega(\mu)} + O(\varepsilon^{1/2}).$$
		For convergent $\int_0^1 \omega(c_1)\cdots\omega(c_n)$, the two regularizing exponentials at the front and the end can be ignored, so
		\begin{equation}\label{aux_10}\int_0^1 \omega(c_1)\cdots\omega(c_n) = \textbf{L}[R^\ast\omega(c_1)\cdots R^\ast\omega(c_n)],\qquad \text{ where }\textbf{L}:= \widetilde{\textbf{I}} e^{B\omega(0)} (g\cdot \widetilde{\textbf{I}}^{-1}).\end{equation}
		It remains to determine the value of $B$. As in the proof of Theorem \ref{graphconnectedtheo}, this can be done by comparing linear terms. Consider
		$$\begin{aligned}\int_0^1 \omega(R(0)) &= \textbf{L}[R^\ast \omega(R(0))] = \textbf{L}[\omega(0)] - \textbf{L}[\omega(\mu^2)] \\ &= B - \widetilde{\textbf{I}}[\omega(\mu^2)] - (g\cdot \widetilde{\textbf{I}}^{-1})[\omega(\mu^2)] \\
			&= B - \int_0^1 \omega(\mu^2) + \int_0^1 \omega(\mu). \end{aligned}$$
		Thus, $$B = \int_0^1 \omega(R(0))+\omega(\mu^2)-\omega(\mu) = -\frac{2\pi i}{5}.$$
		This can also be seen directly since the change of argument on the small circular arc at origin is $-2\pi/5$.
		Equation \ref{aux_10} gives an explicit way to express $\int_0^1 \omega(c_1)\cdots\omega(c_n)$ in terms of multiple polylogarithms at fifth roots of unity. \par
		For instance, consider $\int_0^1 \omega(0)\omega(1)\omega(\alpha)$ with $\alpha = (\sqrt{5}+1)/2$. We have
		$$\begin{aligned}\int_0^1 \omega(0)\omega(1)\omega(\alpha) &= \textbf{L}[R^\ast(\omega(0))R^\ast(\omega(1))R^\ast(\omega(\alpha))] \\&= \textbf{L}[(\omega(\mu)-\omega(\mu^2))(\omega(1)-\omega(\mu^2))(\omega(\mu^4)-\omega(\mu^2))].\end{aligned}$$
		We now split it into $8$ terms and determine each $\textbf{L}[\omega(a_1)\omega(a_2)\omega(a_3)]$. Since each coefficient of $\textbf{L}$ can be effectively expressed as a CMZV, so can our original iterated integral. The computation is quite involved and is best delegated to a computer.
	\end{example}

	\subsection{Mathematica package \textsf{MultipleZetaValues}}\label{package_section1}
	In this subsection, we describe the author's Mathematica package \textsf{MultipleZetaValues}, which implements an effective version of Corollary \ref{main_corollary}. The package can be downloaded at \href{https://www.researchgate.net/publication/357601353}{https://www.researchgate.net/publication/357601353}.
	
	Given a positive integer $N$, let $S = \{0,\infty,1,\mu,\cdots,\mu^{N-1}\}$, where $\mu=e^{2\pi i /N}$, and consider the following collection of subsets of $\mathbb{P}^1$: 
	$$\mathcal{C}_N := \{R(S)| R \text{ is a rational function}, R^{-1}(R(S)) = S \text{ and } \{0,1,\infty\}\subset R(S)\},$$
	By Proposition \ref{S_unit_eq_finite}, $\mathcal{C}_N$ is finite. If $\{a_1,\cdots,a_n\}$ is contained in some element of $\mathcal{C}_N$, then Corollary \ref{main_corollary} gives
	$$\int_0^1 \omega(a_1)\cdots \omega(a_n) \in \CMZV{N}{n}.$$

	\begin{example}
		From Example \ref{level6Ex} above, we saw that $$\int_0^1 \omega(2) \omega(4) = \int_{0<x_2<x_1<1} \frac{dx_1}{x_1-2}\frac{dx_2}{x_2-4}\in \CMZV{6}{2}.$$
		To generate an explicit equality witnessing the containment, we execute the following command in the Mathematica package \textsf{MultipleZetaValues}.
		
		\begin{mmaCell}[functionlocal=y]{Code}
		\mmaDef{MZExpand}[\mmaDef{IterInt}[{2, 4}], "IterIntToCMZV"]
		\end{mmaCell}
		\begin{mmaCell}{Output}
		ColoredMZV[3, \{1, 1\}, \{1, 0\}] - ColoredMZV[3, \{1, 1\}, \{1, 1\}] 
		+ ColoredMZV[3, \{1\}, \{1\}] ColoredMZV[6, \{1\}, \{5\}] - ColoredMZV[6, \{1, 1\}, \{3, 5\}] 
		- ColoredMZV[6, \{1, 1\}, \{5, 4\}] + ColoredMZV[6, \{1, 1\}, \{5, 5\}] 
		- ColoredMZV[3, \{1\}, \{1\}] MultiZeta[\{-1\}] - ColoredMZV[6, \{1\}, \{5\}] MultiZeta[\{-1\}] 
		+ MultiZeta[\{-1\}]^2 + MultiZeta[\{-1, 1\}]
		\end{mmaCell}
		
		Here, \texttt{ColoredMZV[N,\string{s1,...,sn\string},\string{a1,...,an\string}]} is the value of the multiple polylogarithm $\Li_{s_1,\cdots,s_n}(\mu^{a_1},\cdots,\mu^{a_n})$ with $\mu=\exp(2\pi i/N)$. 
		
		For the level $5$ example $\int_0^1 \omega(0) \omega(1) \omega((\sqrt{5}+1)/2)$, one simply executes
		\begin{mmaCell}[functionlocal=y]{Code}
		\mmaDef{MZExpand}[\mmaDef{IterInt}[{0,1,(Sqrt[5]+1)/2}], "IterIntToCMZV"]
		\end{mmaCell}
	\end{example}
	
	\begin{example}\label{IterIntDoableQ_example}
		For positive integers $N\leq 12$, the set $\mathcal{C}_N$ is stored internally in the Mathematica package. For a given set $\{a_1,\cdots,a_n\}$, the command \textsf{IteratedIntDoableQ} checks whether it is contained in some element of $\mathcal{C}_N$ for some $N\leq 12$. For example,
		\begin{mmaCell}[functionlocal=y]{Code}
		{\mmaDef{IterIntDoableQ}[{0,1,2,4}], \mmaDef{IterIntDoableQ}[{(0,1,(Sqrt[5]+ 1)/2, (Sqrt[5]+3)/2}],
		\mmaDef{IterIntDoableQ}[{4Sqrt[5]-8, 5Sqrt[5]-10}]}
		\end{mmaCell}
		\begin{mmaCell}{Output}
		\{6,5,10\}
		\end{mmaCell}
		The output states that $\int_0^1 \omega(a_1)\cdots \omega(a_n) \in \CMZV{6}{n}$ when $a_i \in \{0,1, 2,4\}$; $\int_0^1 \omega(a_1)\cdots \omega(a_n) \in \CMZV{5}{n}$ when $a_i \in \{0,1, (1+\sqrt{5})/2,(3+\sqrt{5})/2\}$; and $\int_0^1 \omega(a_1)\cdots \omega(a_n) \in \CMZV{10}{n}$ when $a_i \in \{0,1, 4\sqrt{5}-8, 5\sqrt{5}-10\}$. \par
		Here are further examples at other levels:
		\begin{mmaCell}[functionlocal=y]{Code}
		{\mmaDef{IterIntDoableQ}[{Csc[Pi/14] Sin[3 Pi/14], -2 Sin[Pi/14]}], 
		\mmaDef{IterIntDoableQ}[{(0,97 + 56 Sqrt[3], 21 + 12 Sqrt[3]}],
		\mmaDef{IterIntDoableQ}[{-1, -I Sqrt[3 + 2 Sqrt[2]], -I Sqrt[3 - 2 Sqrt[2]]}]}
		\end{mmaCell}
		\begin{mmaCell}{Output}
		\{7,12,8\}
		\end{mmaCell}
		
		For each of these examples, one can obtain an explicit CMZV expression for the corresponding iterated integrals. For instance, 
		\begin{mmaCell}[functionlocal=y]{Code}
		\mmaDef{MZExpand}[\mmaDef{IterInt}[{Csc[Pi/14] Sin[3 Pi/14], 0, 1, -2 Sin[Pi/14]}], "IterIntToCMZV"]
		\end{mmaCell}
		gives an explicit reduction of the iterated integral in terms of level 7 CMZVs.
	\end{example}
	
	The algorithm used by the above commands is as follows. For each rational function $R$ in the definition of $\mathcal{C}_N$, one can write down a formal power series $\textbf{L}$ (in general, it has the form appearing in (\ref{aux_9})) such that $\int_0^1 w = \textbf{L}[R^\ast w]$. The program then computes the corresponding coefficient. For any such $R$, the corresponding $\textbf{L}$ is hard-coded into the package.

	\section{$\mathbb{Q}$-relations among CMZVs}
	The three smallest levels for which non-standard relations occur are $N=4, 6$, and $8$. In these cases, the motivic dimension of $\CMZV{N}{n}$ is given by the formula (\ref{deligne_bound}). We describe a class of relations, which we call \textit{$S$-unit relations}, that seems able to produce all non-standard relations for these three levels. \par
	In this section, we use the following notation for differential forms\footnote{The same notation is used by Zhao in \cite{ZhaoStandard, zhao2008multiple} and \cite[Chap.~14]{zhao2016multiple}.}: $$x_0 = \frac{dx}{x} = \omega(0), \quad x_1 = \frac{dx}{1-x}.$$
	Also $$a = \frac{dx}{x} \qquad b_i = \frac{dx}{\mu^{-i}-x} \qquad \mu=e^{2\pi i /N}.$$
	The reason for using two notations for $dx/x$ will soon become clear. Recall the notion of an $N$-unital function defined in the introduction.

	\subsection{$S$-unit relation: examples}\label{nonsd}
	
	We give three illustrative examples. 
	
	\begin{example}\label{S_unit_ex1}
		Let $N=6$ and $\mu=e^{2\pi i/6}$. The two functions $$R_1 = \frac{(1-i\sqrt{3})x}{2(x-\mu^2)^2}, \qquad T_1 = \frac{x^2}{1+x+x^2},$$
		are both $6$-unital, with $$(R_1^\ast x_0, R_1^\ast x_1) = (a+2b_4,b_3-2b_4+b_5),\quad (T_1^\ast x_0, T_1^\ast x_1) = (2a+b_2+b_4, -b_2+b_3-b_4).$$
		Consider the paths $R_1\circ [0,1]$ and $T_1\circ [0,1]$, both of which start at $0$ and end at $1/3$. For $w\in \{x_0,x_1\}^\ast$ not ending in $x_0$ (to ensure convergence), the paths $R_1\circ [0,1]$ and $T_1\circ [0,1]$ are homotopic, and hence
		\begin{figure}[h]
			\centering
			\includegraphics[scale=0.8]{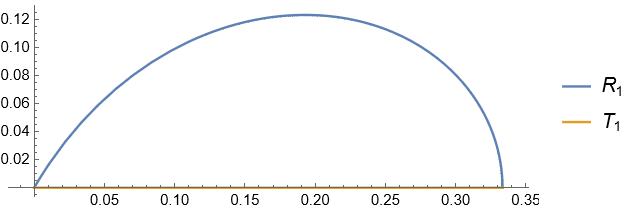}
			\caption{\small The paths $R_1\circ [0,1]$ and $T_1\circ [0,1]$.}
		\end{figure}
		$$\int_{R_1\circ [0,1]} w = \int_{T_1\circ [0,1]} w\iff \int_0^1 R_1^\ast w = \int_0^1 T_1^\ast w.$$
		For example, letting $w = x_0 x_0 x_1$ yields the relation $\int_0^1 u = 0$ with 
		{\small \begin{multline*}u = 4 a a b_2-3 a a b_3+2 a a b_4+a a b_5+2 a b_2 b_2-2 a b_2 b_3+2 a b_2 b_4+2 a b_4 b_2-2 a b_4 b_4+2 a b_4 b_5+2 b_2 a b_2-2 b_2 a b_3\\ +2 b_2 a b_4+2 b_4 a b_2-2 b_4 a b_4+2 b_4 a b_5 +b_2 b_2 b_2-b_2 b_2 b_3+b_2 b_2 b_4+b_2 b_4 b_2-b_2 b_4 b_3+b_2 b_4 b_4\\ +b_4 b_2 b_2-b_4 b_2 b_3+b_4 b_2 b_4+b_4 b_4 b_2+3 b_4 b_4 b_3-7 b_4 b_4 b_4+4 b_4 b_4 b_5.\end{multline*}}
		This is a relation among level-$6$, weight-$3$ CMZVs that cannot be generated by standard relations. 
	\end{example}
	
	\begin{remark}
		We also note that, at present, the only way to determine whether a given relation is non-standard is through explicit computation: one must enumerate all standard relations and then test linear independence against this relation. This involves performing Gaussian elimination on $\mathbb{Q}$-matrices of fairly large size.
	\end{remark}
	
	\begin{example}\label{S_unit_ex2}
		Let $N=6, \mu=e^{2\pi i/6}$. Let $$R_1 = \frac{1-\mu^2 x}{1+x},\qquad T_1 = \frac{2 \left(x-\mu^2\right) \left(x-\mu^5\right)}{\left(1+i \sqrt{3}\right) (x+1)}.$$
		They are both $6$-unital, with $$(R_1^\ast x_0,R_1^\ast x_1) = (b_3-b_2,-a-b_3), \qquad ( T_1^\ast x_0, T_1^\ast x_1) = (-b_1+b_3-b_4,-a-b_3+b_5).$$
		\begin{figure}[h]
			\centering
			\includegraphics[scale=0.8]{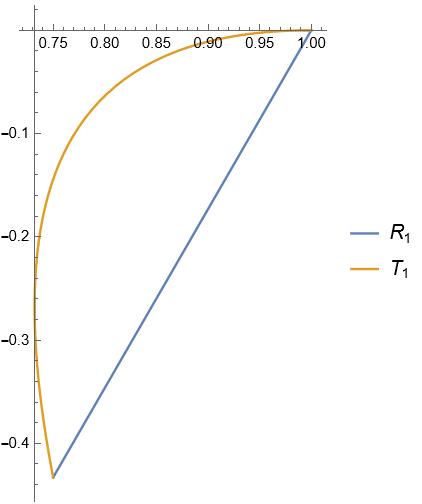}
			\caption{\small The paths $R_1\circ [0,1]$ and $T_1\circ [0,1]$.}
		\end{figure}
		Consider the paths $R_1\circ [0,1]$ and $T_1\circ [0,1]$, both of which start at $1$ and end at $(1-\mu^2)/2$. For $w\in \{x_0,x_1\}^\ast$ not ending in $x_1$ (to ensure convergence), the paths $R_1\circ [0,1]$ and $T_1\circ [0,1]$ are homotopic, and hence
		$$\int_{R_1\circ [0,1]} w = \int_{T_1\circ [0,1]} w\iff \int_0^1 R_1^\ast \omega = \int_0^1 T_1^\ast \omega.$$
		Taking $w$ to be certain words of weight $3$ yields new non-standard relations that are independent of those in the previous example.
	\end{example}
	
	In the above two examples, we considered only words for which both integrals converge (the ``finite'' version), but it is not difficult to regularize them. Once the level $N$ is clear from the context, recall the group-like element $\widetilde{\textbf{I}}$ defined at the start of Subsection \ref{expcomp}. For an $N$-unital function $R$, define
	$$\textbf{I}_R := \sum_{w\in \{x_0,x_1\}^\ast} \widetilde{\textbf{I}}[R^\ast w]w.$$
	Then for $(R_1,T_1)$ in Example \ref{S_unit_ex1}, we have
	$$\textbf{I}_{R_1} = \textbf{I}_{T_1} e^{Ax_0};$$
	for $(R_1,T_1)$ in Example \ref{S_unit_ex2}, we have
	$$\textbf{I}_{R_1} = \textbf{I}_{T_1} e^{Bx_1}.$$
	The constants $A, B\in \mathbb{C}$ in both examples can be found by comparing the weight-$1$ coefficients. The above two displayed equations are the regularized relations from Examples \ref{S_unit_ex1} and \ref{S_unit_ex2}. In some cases, regularization is unavoidable, as the next example shows. 	
	\begin{example}
		Let $N=8$ and $\mu=e^{2\pi i/8}$. Let $R_1,R_2,T_1,T_2$ be $8$-unital functions such that\footnote{The rational functions $R_i, T_i$ are uniquely determined by these data.}
		$$R_1^\ast(x_0,x_1) = (-a-b_6,a+b_5)\qquad R_2^\ast(x_0,x_1) = (b_2+b_3-2b_5,-b_2-b_3+b_4+b_7)$$
		$$T_2^\ast(x_0,x_1) = (-a-b_6,a+b_7)\qquad T_1^\ast(x_0,x_1) = (-b_5+b_6,b_4-b_6)$$
		Then $R_1(0)=T_1(0) = \infty$ and $R_2(1) = T_2(1)$. Consider the two composite paths: $R_1\circ [0,1]$ followed by $R_2\circ [0,1]$, and $T_1\circ [0,1]$ followed by $T_2\circ [0,1]$. They both start at $\infty$ and end at the same point, which is neither $0$ nor $1$.
		\begin{figure}[H]
			\centering
			\includegraphics[scale=0.9]{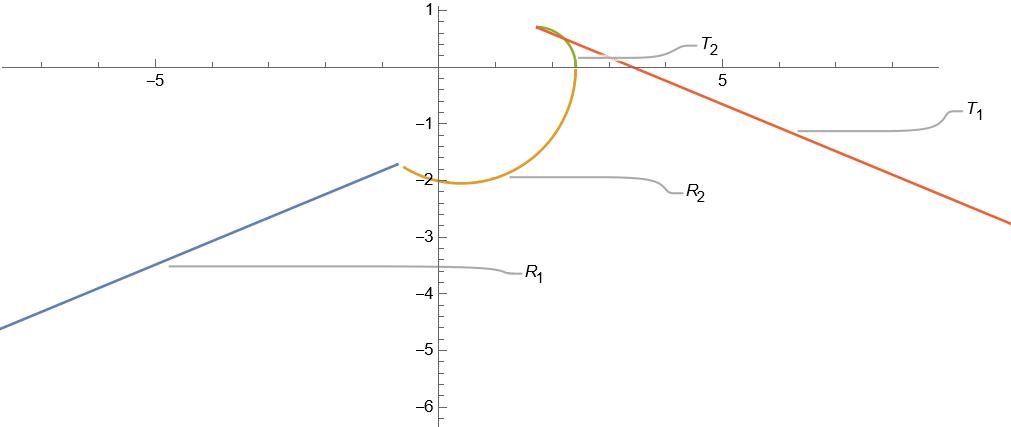}
			\caption{\small The paths $R_i\circ [0,1]$ and $T_i\circ [0,1]$.}
		\end{figure}
	\end{example}
	Lemma \ref{radiusinf} gives the regularizing exponential at $\infty$, and we have \begin{equation}\label{temp1}\textbf{I}_{R_2}\textbf{I}_{R_1} e^{A(x_0-x_1)} = \textbf{I}_{T_2}\textbf{I}_{T_1}\end{equation}
	for some constant $A$. Comparing linear terms, one sees that $A = -3\pi i/4$. Comparing the coefficients on both sides of (\ref{temp1}) gives non-standard relations among level-$8$ CMZVs.
	
	\subsection{$S$-unit relations: general formulation}
	Let $R_1,R_2,\cdots,R_r$ be $N$-unital functions such that
	\begin{itemize}
		\item $a_0 := R_1(0) \in \{0,1,\infty\}$;
		\item $a_i := R_i(1) = R_{i+1}(0)$ for $i=1,\cdots,r-1$; and
		\item $a_r := R_r(1) = R_1(0)$.
	\end{itemize}
	That is, the paths $R_i\circ [0,1]$, $i=1,\cdots,r$, can be composed to form a loop that starts and ends at $a_0 \in \{0,1,\infty\}$. Let $$X(a) := \begin{cases}
		x_0 & \text{ if }a=0, \\
		x_1 & \text{ if } a=1, \\
		x_0 - x_1 & \text{ if } a=\infty,\\
		0& \text{ otherwise}.
	\end{cases}$$
	Then we have \begin{equation}\label{aux_11}\textbf{I}_{R_r} e^{A_{r-1} X(a_{r-1})}\cdots \textbf{I}_{R_2} e^{A_1 X(a_1)} \textbf{I}_{R_1} e^{A_0 X(a_0)}\widetilde{\textbf{I}}_{\rho} = 1,\end{equation}
	where $\rho$ is a loop based at $a_0$ and $\widetilde{\textbf{I}}_{\rho}$ is defined in Proposition \ref{mainregtheorem}. In the above three examples, $\rho$ is null-homotopic, so this term can be ignored. The constants $A_i \in \CMZV{N}{1}$ can be found by comparing weight-$1$ terms. Equation (\ref{aux_11}) is what we mean by an $S$-unit relation in full generality. 
	
	Let $$\overline{\CMZV{N}{n}} := \sum_{1\leq k <n} \CMZV{N}{k}\CMZV{N}{n-k}.$$ By Theorem \ref{homotopymodzero}, the weight-$n$ coefficient of $\widetilde{\textbf{I}}_{\rho}$ is in $\overline{\CMZV{N}{n}}$. It is more elegant to write (\ref{aux_11}) modulo $\overline{\CMZV{N}{n}}$: for any word $w\in \{x_0,x_1\}^\ast$ of weight $n$,
	$$\widetilde{\textbf{I}}[R_r^\ast w + \cdots + R_1^\ast w] \equiv 0 \pmod{\overline{\CMZV{N}{n}}}.$$
	
	By combining standard relations with $S$-unit relations, Deligne's bound can be reached for the following levels and weights:\footnote{Non-standard relations exist for these four levels.}
	\begin{itemize}
		\item Level $6$, weight $\leq 5$;
		\item Level $8$, weight $\leq 4$;
		\item Level $10$, weight $\leq 3$;
		\item Level $12$, weight $\leq 3$.
	\end{itemize}
	However, these relations seem unable to reach Deligne's bound for level $10$, weight $4$.\footnote{In this case, there are $72$ non-standard relations. The $S$-unit relations give $70$ of them, while the remaining $2$ remain elusive.} Moreover, for level $9$, weight $3$, for which there are $3$ non-standard relations, the $S$-unit relations do not produce anything new.
	\begin{conjecture}
		For levels $N=6,8$, all non-standard relations come from $S$-unit relations.
	\end{conjecture}

	\subsection{CMZV database}
	The relationship between Deligne's bound (Theorem \ref{deligne_bound}) and non-standard relations is quite complex for general $N$. We summarize it as follows (see \cite{zhao2008multiple, ZhaoStandard} for more details):
	\begin{itemize}
		\item When $N=1,2,3,4,6,8$, Deligne's bound is tight. If, moreover, $N=4,6,8$, non-standard relations exist.
		\item When $N=p^n$, where $n\geq 1$ and $p\geq 5$ is prime, Deligne's bound is known not to be tight, while non-standard relations do not seem to exist.
		\item For other $N$, non-standard relations seem to exist and we do not know whether Deligne's bound is tight.
	\end{itemize}
	
	The Mathematica package \textsf{MultipleZetaValues} includes a database of CMZVs of small weights $n$ and levels $N$. In the current version (version 1.2.0), the available pairs $(N,n)$ are:
	\begin{itemize}
		\item $n\leq 14$ for level 1;
		\item $n\leq 8$ for level 2;
		\item $n\leq 5$ for level 3;
		\item $n\leq 6$ for level 4;
		\item $n\leq 4$ for level 5;
		\item $n\leq 5$ for level 6;
		\item $n\leq 4$ for level 8;
		\item $n\leq 3$ for levels 7, 10, and 12.
	\end{itemize}
	For each $(N,n)$ mentioned above, the package \textsf{MultipleZetaValues} contains an explicit list of complex numbers $\mathcal{B}^{N}_n$ such that $\CMZV{N}{n}$ equals the $\mathbb{Q}$-span of $\mathcal{B}^{N}_n$. These numbers constitute a basis assuming Grothendieck's period conjecture.
	
	\begin{example}
		To display the explicit constants in $\mathcal{B}^N_n$, one simply executes \textsf{MZBasis[N,n]}. 
		\begin{mmaCell}[functionlocal=y]{Code}
		\mmaDef{MZBasis}[6,2]
		\end{mmaCell}
		\begin{mmaCell}{Output}
		\{I Sqrt[3] DirichletL[3, 2, 2], PolyLog[2, 1/4], Pi^2, I Pi Log[3], I Pi Log[2], 
		Log[3]^2, Log[2] Log[3], Log[2]^2\}
		\end{mmaCell}
		This gives an explicit basis for the level-$6$, weight-$2$ CMZVs, where \textsf{DirichletL[3, 2, s]} represents the Dirichlet $L$-function $L_{-3}(s) = \sum_{n\geq 0} (\frac{1}{(3n+1)^s} - \frac{1}{(3n+2)^s})$. 
		
		The following expresses multiple polylogarithms at roots of unity in terms of $\mathcal{B}^N_n$. 
		\begin{mmaCell}[functionlocal=y]{Code}
		{\mmaDef{ColoredMZV}[2,{1,1,1},{1,1,0}], \mmaDef{ColoredMZV}[2,{2,1,1},{0,0,1}],
		\mmaDef{ColoredMZV}[3,{1,1,1},{2,1,1}]} // \mmaDef{MZExpand}
	\end{mmaCell}
		which gives 	$$\begin{aligned}
		\Li_{1,1,1}(-1,-1,1) &= -\frac{7 \zeta (3)}{8}-\frac{1}{6} \log ^3(2)+\frac{1}{12} \pi ^2 \log (2) \\
		\Li_{2,1,1}(1,1,-1) &= -\text{Li}_4\left(\frac{1}{2}\right)-\frac{7}{8} \zeta (3) \log (2)+\frac{\pi ^4}{80}-\frac{1}{24} \log ^4(2)-\frac{1}{12} \pi ^2 \log ^2(2) \\
		\Li_{1,1,1}(\mu^2,\mu,\mu) &= \frac{\pi  L_{-3}(2)}{2 \sqrt{3}}+\frac{i\sqrt{3}}{4} L_{-3}(2) \log (3)-\frac{2 \zeta (3)}{3}-\frac{5 i \pi ^3}{432}-\frac{1}{48} \log ^3(3)-\frac{1}{48} i \pi  \log ^2(3)+\frac{5}{144} \pi ^2 \log (3)
	\end{aligned}$$
\end{example}

\begin{example}
		By Subsection \ref{package_section1}, $\int_0^1 \omega(0) \omega(2) \omega(4) \in \CMZV{6}{3}$. To express it explicitly using constants in $\mathcal{B}^6_3$, we execute
	\begin{mmaCell}[functionlocal=y]{Code}
		\mmaDef{IterInt}[{0, 2, 4}] // \mmaDef{MZExpand}
	\end{mmaCell}
	\begin{mmaCell}{Output}
		1/12 Pi^2 Log[2] - Log[2]^3/3 - 1/4 PolyLog[3, 1/4] - (7 Zeta[3])/24
	\end{mmaCell}
	
	A slightly less straightforward example is 
	\begin{mmaCell}[functionlocal=y]{Code}
		\mmaDef{IterInt}[{0, 1, (3 + Sqrt[5])/2, 1}] // \mmaDef{MZExpand}
	\end{mmaCell}
	\begin{mmaCell}{Output}
		-((11 Pi^4)/450) + 1/5 Pi^2 Log[GoldenRatio]^2 - Log[GoldenRatio]^4/8 
		- 3/8 PolyLog[4, 1/2 (3 - Sqrt[5])] + 3 PolyLog[4, 1/2 (-1 + Sqrt[5])]
	\end{mmaCell}
\end{example}

We try to make $\mathcal{B}^{N}_{n}$ consist of ``elementary constants.'' This means that priority is given to constants such as $$\log \alpha,\quad \zeta(n),\quad L(\chi,n),\quad \Li_n(\alpha), \qquad \alpha\in \overline{\mathbb{Q}}.$$ For example,
$$\begin{aligned}
	\mathcal{B}^2_1 &= \left\{\log 2\right\} \qquad \mathcal{B}^2_2 = \left\{\pi ^2, \log^2 (2)\right\} \qquad \mathcal{B}^2_3= \left\{\zeta (3),\pi ^2 \log (2),\log ^3(2)\right\} \\
	\mathcal{B}^2_4 &= \left\{\text{Li}_4\left(\frac{1}{2}\right),\zeta (3) \log (2),\pi ^4,\pi ^2 \log ^2(2),\log ^4(2)\right\} \\
	\mathcal{B}^3_1 &= \{i\pi,\log 3\} \qquad \mathcal{B}_{3,2}= \left\{i \sqrt{3} L_{-3}(2),\pi ^2,i \pi  \log (3),\log ^2(3)\right\} \\
	\mathcal{B}^3_3 &= \left\{\zeta (3),i \Im\left(\text{Li}_3\left(\frac{1}{2}+\frac{i}{2 \sqrt{3}}\right)\right),\sqrt{3} \pi  L_{-3}(2),i \sqrt{3} L_{-3}(2) \log (3),i \pi ^3,\pi ^2 \log (3),i \pi  \log ^2(3),\log ^3(3)\right\} \\
	\mathcal{B}^4_3 &= \left\{\zeta (3),i \Im\left(\text{Li}_3\left(\frac{1}{2}+\frac{i}{2}\right)\right),\pi  G,i G \log (2),i \pi ^3,\pi ^2 \log (2),i \pi  \log ^2(2),\log ^3(2)\right\}
\end{aligned}$$
where $G = L_{-4}(2)$ is Catalan's constant. For large $N$ and $n$, such a simple basis may not be available; in that case, we randomly choose some constants of higher depth. One motivation for favoring elementary constants is that they allow us to prove certain classical identities quickly (see the next two sections). 

As mentioned above, the Mathematica package \textsf{MultipleZetaValues} can be regarded as a CMZV database of low levels and weights. We compare it with other databases in the literature.

\begin{itemize}[leftmargin=*]
	\item The MZV Datamine \cite{MZVdatamine}: this is the earliest and still the most comprehensive database for level-$1$ and level-$2$ CMZVs. 
	\item Ablinger \cite{ablinger2011harmonic, ablinger2014iterated} wrote a Mathematica package \textsf{HarmonicSums}, which focuses on iterated integrals whose differential forms are $dx/\Phi_N(x)$, with $\Phi_N(x)$ the $N$-th cyclotomic polynomial. They form a subspace of the level-$N$ CMZVs. The package has a database for $N=1,2,4,6$. Ablinger already noted that his relations are incomplete\footnote{That is, numerical relations exist that are not derivable by his methodology.}. This article supplies the missing relations by finding all relations in the larger space $\CMZV{6}{}$. Ablinger's package also has functionality similar to \textsf{MZIntegrate}, which we will use in the last section.
	\item Duhr and Dulat \cite{duhr2019polylogtools} wrote a Mathematica package \textsf{PolyLogTools}, focusing on the coalgebra structure of iterated integrals. It also has a limited database of special values that covers some iterated integrals of levels $4$ and $6$.
	\item Maître \cite{maitre2006hpl, maitre2012extension} developed a Mathematica package \textsf{HPL} for iterated integrals, which is required by the package of Duhr and Dulat. 
	\item Smirnov, Smirnov, and Henn \cite{henn2017evaluating} compiled an \textit{empirical} database for CMZVs of level $6$ and weight $\leq 6$.
	\item Panzer \cite{panzer2015algorithms} wrote a Maple package \textsf{HyperInt} for generalized polylogarithms. It emphasizes algebraic manipulations and is not intended to provide a database function.
\end{itemize}

\subsection{$\MZV{S}{n}$ for other $S$}\label{generalS}
Three cases for $S$ not equivalent to $\{0,\infty,1,\mu,\cdots,\mu^{N-1}\}$ have been investigated empirically. 
\begin{itemize}[leftmargin=*]
	\item A \textit{multiple Deligne value (MDV)} \cite{broadhurst2014multiple} belongs to the $\mathbb{Q}$-space spanned by convergent integrals of the form
	$$\int_0^1 \omega(a_1)\cdots\omega(a_n),\qquad a_i\in \{0,1,e^{2\pi i /6}\}.$$
	A related class of objects, known as \textit{multiple Clausen values} \cite{borwein2001central}, is defined to be the $\mathbb{Q}$-span of $\Li_{s_1,\cdots,s_k}(e^{\pm i\pi/3})$. It can be shown that the space of multiple Clausen values is the same as the space of MDVs \cite[Section~2.4]{panzer2015algorithms}. Let $S = \{0,1,\infty,e^{2\pi i /6}\}$. These points are the vertices of a regular tetrahedron with symmetry group $G = A_4$, so, by Example \ref{complete_edge_ex}, $\{(0,1)\}$ is a set of complete edges. Theorem \ref{graphconnectedtheo} implies that $\MZV{S}{n}$ coincides with the space of MDVs. Broadhurst conjectured that its dimensions are given by
	$$\sum_{n\geq 0} (\dim_\mathbb{Q} \textsf{MDV}_n)t^n  \stackrel{?}{=} \frac{1}{1-t-t^2}.$$
	From Example \ref{level6Ex_Deligne_value}, we know that $\MZV{S}{n} \subset \CMZV{3}{n}$. 
	\item A \textit{multiple Landen value (MLV)} \cite{broadhurst2015multiple} belongs to the $\mathbb{Q}$-space spanned by convergent integrals of the form
	$$\int_0^1 \omega(a_1)\cdots\omega(a_n),\qquad a_i\in \left\{0,1,\frac{1+\sqrt{5}}{2},\frac{3+\sqrt{5}}{2}\right\}.$$
	Let $S = \{0,\infty,1,\frac{1+\sqrt{5}}{2},\frac{3+\sqrt{5}}{2}\}$. The symmetry group\footnote{In fact, up to a Möbius transformation, $S$ is a planar regular pentagon.} of $S$ is dihedral of order $10$, and $\{(0,1)\}$ is a set of complete edges. Theorem \ref{graphconnectedtheo} implies that $\MZV{S}{n}$ coincides with the space of MLVs\footnote{Modulo constants of weight $1$.}. Broadhurst conjectured that
	$$\sum_{n\geq 0} (\dim_\mathbb{Q} \textsf{MLV}_n)t^n  \stackrel{?}{=} \frac{1}{1-t-t^2-t^3}.$$
	From Example \ref{level5Ex}, we know that $\MZV{S}{n} \subset \CMZV{5}{n}$. 
\end{itemize}

Broadhurst raised the problem of rigorously determining the $\mathbb{Q}$-linear relations among MDVs and MLVs. From their definitions, the only immediately available relations are the shuffle relations, which are far from sufficient to account for the conjectural dimensions of these spaces. However, as we have seen, both MDVs and MLVs embed into the spaces of CMZVs of levels $3$ and $5$, respectively. Since linear relations in these CMZV spaces are understood\footnote{In particular, there should be no non-standard relations at these levels.}, all relations between MDVs and MLVs can now be found in these larger spaces. It would nevertheless be interesting to see whether such relations can be obtained directly, without passing to an ambient CMZV space. Recently, this has been done for MDVs via their tetrahedral symmetry \cite[Section~3]{sun2026multiple}.

Broadhurst also investigated the so-called \textit{multiple Watson values} \cite{broadhurst2015tests}, which are convergent iterated integrals of the form
$$\int_0^1 \omega(a_1)\cdots\omega(a_n),\qquad a_i\in \left\{0,1,\gamma,\gamma^2,\frac{\gamma}{1+\gamma},\frac{\gamma^2}{1-\gamma}\right\},\qquad \gamma = 2\sin(\frac{\pi}{14}).$$
Our theory does not apply here because the corresponding $S$ has a trivial symmetry group, although it might be related to level-$14$ CMZVs. \par

In the next section, we briefly investigate another case in which $S$ is the set of $12$ vertices of a regular icosahedron. 

% The \textit{shuffle relation} generalizes verbatim to $\MZV{S}{n}$. On the other hand, the parallel of \textit{stuffle relation} is not immediate; in fact, in case of general $S$, for each cyclic subgroup of $G$, up to conjugation, there is a version of stuffle algebra, form which stuffle relations can be deduced. The usual stuffle algebra of CMZV corresponding to the cyclic subgroup $\{1,\mu,\cdots,\mu^{N-1}\}$ of $G$, with $G$ being dihedral of order $2N$. More about this might appear in a later article. 

\section{Polylogarithm identities}
\subsection{Special values of multiple polylogarithms}

When we know all relations among CMZVs of a given weight and level, every purported identity in this space can be checked. 

\begin{example}
	The following three famous closed-form evaluations of dilogarithms and trilogarithms
	$$\text{Li}_2\left(\frac{\sqrt{5}-1}{2}\right)=\frac{\pi ^2}{10}-\log ^2(\phi )\qquad \text{Li}_2\left(\frac{3-\sqrt{5}}{2}\right)=\frac{\pi ^2}{15}-\log ^2(\phi )$$
	$$\text{Li}_3\left(\frac{3-\sqrt{5}}{2}\right)=\frac{4 \zeta (3)}{5}+\frac{2 \log ^3(\phi )}{3}-\frac{2}{15} \pi ^2 \log (\phi )$$
	can now be checked easily because both sides are level-$5$ CMZVs.
\end{example}

\begin{example}
	We also have an analogue for multiple polylogarithms, where $\rho = (\sqrt{5}-1)/2$: 
	$$\Li_{1,1,1,1}(\rho,1,1,\rho) = 2 \text{Li}_3(\rho) \log (\phi )-\frac{1}{4} \text{Li}_4(\rho^2)+2 \text{Li}_4(\rho)-\frac{2}{5} \zeta (3) \log (\phi )+\frac{5 \log ^4(\phi )}{8}+\frac{7}{60} \pi ^2 \log ^2(\phi )-\frac{7 \pi ^4}{360}$$
	$$\Li_{2,2}\left(\frac{1}{2},2\right) = \sum_{n=1}^\infty \frac{1}{n^22^n}\sum_{m=1}^{n-1}\frac{2^m}{m^2} = -3 \text{Li}_4\left(\frac{1}{2}\right)+\frac{7 \pi ^4}{288}-\frac{1}{8} \log ^4(2)-\frac{1}{8} \pi ^2 \log ^2(2)$$
\end{example}

\begin{example}\label{coxeterladder}
	The following three dilogarithm ladders due to Coxeter \cite{coxeter1935functions} are classical:
	$$\begin{aligned}
		\text{Li}_2\left(\rho ^6\right)&=4 \text{Li}_2\left(\rho ^3\right)+3 \text{Li}_2\left(\rho ^2\right)-6 \text{Li}_2(\rho )+\frac{7 \pi ^2}{30} \\
		\text{Li}_2\left(\rho ^{12}\right)&=2 \text{Li}_2\left(\rho ^6\right)+3 \text{Li}_2\left(\rho ^4\right)+4 \text{Li}_2\left(\rho ^3\right)-6 \text{Li}_2\left(\rho ^2\right)+\frac{\pi ^2}{10} \\
		\text{Li}_2\left(\rho ^{20}\right)&=2 \text{Li}_2\left(\rho ^{10}\right)+15 \text{Li}_2\left(\rho ^4\right)-10 \text{Li}_2\left(\rho ^2\right)+\frac{\pi ^2}{5} 
	\end{aligned}$$
	with $\rho = \phi^{-1} = (\sqrt{5}-1)/2$. Both sides of the first identity are level-$10$ CMZVs, so the identity can now be verified routinely. We will prove the last one later in Section \ref{Coxeter3rd}. The middle one remains elusive from our perspective. We also have the following ladders, in which both sides are level-$10$ CMZVs. 
	$$\text{Li}_3(\rho^6)-8 \text{Li}_3(\rho^3)-6 \text{Li}_3(\rho) = \frac{3 \zeta (3)}{5}-4 \log ^3(\phi )+\frac{2}{5} \pi ^2 \log (\phi )$$
	
	$$36 \text{Li}_4(\rho)-\frac{9 \text{Li}_4(\rho^2)}{4}-16 \text{Li}_4(\rho^3)+\text{Li}_4(\rho^6) = \frac{27 \log ^4(\phi )}{4}-\frac{3}{2} \pi ^2 \log ^2(\phi )+\frac{2 \pi ^4}{9}$$
	The above two identities were hinted at in \cite[p.~44]{lewin1991structural} via classical ladder techniques.
\end{example}

\begin{example}
	The following dilogarithm identity was discovered by Watson in 1937 \cite{watson1937note}:
	$$\Li_2(\alpha) - \Li_2(\alpha^2) = \frac{\pi^2}{42}+\log^2 \alpha, \qquad \alpha = \frac{1}{2}\sec \frac{2\pi}{7}.$$
	This can also be handled by our approach because both sides are level-$7$ CMZVs.
\end{example}

\begin{example}
	The following ladders involving powers of $-1/2$ are amenable to our approach; they have level $6$ and weights $4$ and $5$, respectively:
	$$\Li_4\left(\frac{-1}{8}\right) = -12 \text{Li}_4\left(\frac{1}{2}\right)+\frac{27 \text{Li}_4\left(\frac{1}{4}\right)}{4}+\frac{\pi ^4}{18}+\frac{5 \log ^4(2)}{8}-\frac{1}{4} \pi ^2 \log ^2(2)$$
	$$\Li_5\left(\frac{-1}{8}\right) = -36 \text{Li}_5\left(\frac{1}{2}\right)+\frac{81 \text{Li}_5\left(\frac{1}{4}\right)}{8}+\frac{403 \zeta (5)}{16}-\frac{3}{8} \log ^5(2)+\frac{1}{4} \pi ^2 \log ^3(2)-\frac{1}{6} \pi ^4 \log (2)$$
	
	The corresponding generalization to weight $6$ was mentioned in \cite{charlton2019}:
	\begin{multline}\label{level2conj}\zeta(5,\bar{1}) = \frac{36 \text{Li}_6\left(\frac{1}{2}\right)}{13}-\frac{81 \text{Li}_6\left(\frac{1}{4}\right)}{208}+\frac{\text{Li}_6\left(-\frac{1}{8}\right)}{39}\\ +\frac{3 \zeta (3)^2}{8}+\frac{31}{16} \zeta (5) \log (2)-\frac{1787 \pi ^6}{589680}-\frac{1}{208} \log ^6(2)+\frac{1}{208} \pi ^2 \log ^4(2)-\frac{1}{156} \pi ^4 \log ^2(2)\end{multline}
	Here, $\zeta(5,\bar{1}) = \sum_{i>j\geq 1} \frac{(-1)^j}{i^5j} \in \CMZV{2}{6}$. A closely related ladder, but in can be found in \cite{broadhurst1998polylogarithmic}. The above three displayed equations verify Conjecture \ref{galois_conj}(d) up to $n= 6$. 
\end{example}

\begin{example}
	We show that, for any rational number $a/b$, $$\Li_{s_1,\cdots,s_n}\left(\frac{a}{b}\right) \in \CMZV{N}{s_1+\cdots+s_n},$$
	where the level is $N = \text{lcm}(a,b,a-b)$. To see this, let $R(x) = \frac{1-x^a}{1-x^b}$. Then $R^{-1}(0)$ is a subset of the $a$-th roots of unity, $R^{-1}(1)$ is a subset of the $|a-b|$-th roots of unity, and $R^{-1}(\infty)$ is a subset of the $b$-th roots of unity. Let $x_0 = dx/x$ and $x_1 = dx/(1-x)$, and let $\omega$ be a word in $x_0$ and $x_1$. Then
	$$\int_{R\circ [0,1]} \omega = \int_0^1 R^\ast \omega \in \CMZV{N}{}.$$
	Since $R(0) = 1$ and $R(1) = a/b$, we conclude that $\Li_{s_1,\cdots,s_n}(\frac{a}{b})$ is in the indicated space.
\end{example}

\begin{remark}
	(a) The value of $N$ in the above example is not optimal. For example, $a/b = 8/9$ gives $N=72$, but Table \ref{polylog_table} shows that $6$ is already sufficient. \\
	(b) Using the method of the recent work \cite{charlton2025multiple}, it seems that one could generalize the above conclusion to multiple polylogarithms: for $r_i \in \mathbb{Q}$, $\Li_{s_1,\cdots,s_n}(r_1,\cdots,r_n) \in \CMZV{N}{s_1+\cdots+s_n}$ for some $N$.
\end{remark}

% For level $N$ CMZV, the possible $z$'s in \ref{zinpolylog} all satisfy the following: both $z$ and $1-z$ are of the form $$\mu^i \prod_j (1-\mu^{a_j})^{b_j}$$ for integers $i,a_i,b_i$, which can be readily seen by looking at $\Li_1(z) = -\log(1-z)$ and $\Li_1(1-z) = -\log z$. In particular, for such $z$, both $z$ and $1-z$ are $S$-units in cyclotomic field $\mathbb{Q}(e^{2\pi i /N})$, with\footnote{not to be confused with the previous common notation that $S$ is a finite subset of $\mathbb{P}^1$} $S$ set of prime ideals lying over prime divisors of $N$. These are solutions to the $S$-\textit{unit equation}, which there are only finitely many solutions \cite{siegel2014einige}, and when $N$ is small, these two set seems to coincide. When $N=10$, there exist $S$-units which seems not to be deducible from \ref{zinpolylog}: $-\rho^{10}, -\rho^6$, we will explain the first one when we prove Coxeter's third ladder in \ref{Coxeter3rd}.

\subsection{Icosahedral MZVs and Coxeter's ladder}
Let $\mathcal{I}$ be the set of $12$ vertices of a regular icosahedron embedded in the Riemann sphere $\mathbb{P}^1$. We investigate the space $\MZV{\mathcal{I}}{}$ and prove Coxeter's third ladder. Explicitly, $$\mathcal{I} = \{0,\infty,\rho\mu^i,-\rho^{-1}\mu^i | 0\leq i \leq 4\}, \qquad \rho = \frac{\sqrt{5}-1}{2}, \qquad \mu=e^{2\pi i /5}.$$
Note that $1\notin \mathcal{I}$. Let $G\cong A_5$ be the group of Möbius transformations that permute $\mathcal{I}$. 

\begin{lemma} We have $\Li_{s_1,\cdots,s_n}(-\rho^{10}) \in \MZV{\mathcal{I}}{s_1+\cdots+s_n}$.
\end{lemma}
\begin{proof}
	Let $R(z)=z^5$. Since $\mathcal{I}$ is invariant under multiplication by $e^{2\pi i /5}$, we have $R^{-1}(R(\mathcal{I})) \subset \mathcal{I}$. Proposition \ref{generalStrans} then implies that, for $R(\mathcal{I}) = \{0,\infty,\rho^5,-\rho^{-5}\}$, we have $\MZV{R(\mathcal{I})}{n}\subset \MZV{\mathcal{I}}{n}$.
	After dividing by $-\rho^{-5}$ (again using Möbius invariance), we can replace $R(\mathcal{I})$ by $S:=\{0,\infty,1,-\rho^{10}\}$; we still have $\MZV{S}{n}\subset \MZV{\mathcal{I}}{n}$. For any complex number $z$, $\MZV{\{0,1,\infty,z\}}{}$ contains all generalized polylogarithms at $z$, which completes the proof.
\end{proof}

\begin{lemma} The weight-$1$ space $\MZV{\mathcal{I}}{1}$ has the $\mathbb{Q}$-basis $\{2\pi i, \log 5, \log \rho\}$.
\end{lemma}
\begin{proof}
	This follows by computing the cross-ratios of $4$-tuples in $\mathcal{I}$ and applying Lemma \ref{level1span}.
\end{proof}

\begin{theorem}
	The spaces of icosahedral MZVs of weights $1$ and $2$ are subspaces of the corresponding spaces of level-$10$ CMZVs; that is,
	$$\MZV{\mathcal{I}}{1} \subset \CMZV{10}{1},\qquad \MZV{\mathcal{I}}{2} \subset \CMZV{10}{2}.$$
\end{theorem}
\begin{proof}
	The weight-$1$ case is straightforward: $\MZV{\mathcal{I}}{1}$ is spanned by $\{2\pi i ,\log 5,\log \rho\}$, while $\CMZV{10}{1}$ is spanned by $\{2\pi i ,\log 5,\log \rho,\log 2\}$. The weight-$2$ inclusion requires more work. Set $S = \{0,\infty,1,e^{2\pi i /10},\cdots,e^{2\pi i \times 9/10}\}$. We claim that, for any $a_i,a_j \in \mathcal{I}/\rho$, there exists a rational function $R_{ij}$ such that
	\begin{equation}\label{aux_1}R_{ij}^{-1}(R_{ij}(S)) = S,\qquad \{0,\infty,1,a_i,a_j\} \subset R_{ij}(S).\end{equation}
	Proposition \ref{S_unit_eq_finite} implies that only finitely many computations are needed to verify the above assertion. Then $$\int_0^1 \omega(a_i)\omega(a_j) \in \MZV{R_{ij}(S)}{2} \subset \MZV{S}{2} = \CMZV{10}{2},$$
	by Proposition \ref{generalStrans}. Finally, for any two distinct elements $t,u \in \mathcal{I}$, the geometric interpretation of the action of $G$ on $\mathcal{I}$ makes it clear that $\{(t,u)\}$ is a set of complete edges. Therefore, by Theorem \ref{graphconnectedtheo}, the $\mathbb{Q}$-span of $\int_0^1 \omega(a_i)\omega(a_j)$, where $a_i, a_j$ range over all elements of $\mathcal{I}/\rho$ for which the integral converges, is $\MZV{\mathcal{I}}{2}$. This completes the proof. 
\end{proof}

The exhaustive verification in the above proof can be delegated to the Mathematica package function \textsf{IterIntDoableQ} (see Example \ref{IterIntDoableQ_example}) by executing the following code:
\begin{mmaCell}[functionlocal=list]{Code}
		Block[{list}, list = Join[{0}, 
		GoldenRatio^(-1)*Table[Exp[2 Pi*I*i/5], {i, 0, 4}], 
		-GoldenRatio* Table[Exp[2 Pi*I*i/5], {i, 0, 4}]]; list = list/list[[-1]]; 
		\mmaDef{IterIntDoableQ} /@ Subsets[list, {2}]]
\end{mmaCell}
The output of these commands consists entirely of integers, which proves our assertion. The proof fails for weight $n\geq 3$ because the statement corresponding to (\ref{aux_1}) is false in weight $3$; the relationship between $\MZV{\mathcal{I}}{n}$ and $\CMZV{10}{n}$ is not known.

\begin{corollary}\label{Coxeter3rd}The following identity holds:
	$$\text{Li}_2\left(\rho ^{20}\right)-2 \text{Li}_2\left(\rho ^{10}\right)=15 \text{Li}_2\left(\rho ^4\right)-10 \text{Li}_2\left(\rho ^2\right)+\frac{\pi ^2}{5}.$$
\end{corollary}
\begin{proof}
	Note that $\Li_2(\rho^{20}) - 2\Li_2(\rho^{10}) = 2\Li_2(-\rho^{10})$. Therefore, the left-hand side is in $\MZV{\mathcal{I}}{2}$ by our first lemma, and it is also in $\CMZV{10}{2}$ by the preceding theorem. Thus, the above identity is an equality in $\CMZV{10}{2}$ and can be checked effectively.
\end{proof}

\mycomment{
	Given the weight $2$ subspace inclusion in previous theorem, it would be interesting to ask
	\begin{problem}
		What is the relationship between $\MZV{\mathcal{I}}{n}$ and $\CMZV{10}{n}$? Is the former one subspace of the latter?
	\end{problem}

	\begin{problem}
		Are $\Li_{s_1,\cdots,s_n}(-\rho^6), \Li_{s_1,\cdots,s_n}(\rho^{12}),\Li_{s_1,\cdots,s_n}(-\rho^{10}) \in \CMZV{10}{}$? 
	\end{problem}
}

\section{Applications to Apéry-type infinite series}
Here, we convert some series into iterated integrals and then into CMZVs. When they have weights and levels for which all $\mathbb{Q}$-relations are known, we obtain closed-form evaluations of the series. 

We first record the integration kernel $$\int_0^1 x^{n-1}(1-x)^n dx = \frac{1}{n}\binom{2n}{n}^{-1}.$$

\begin{proposition}
	For $n\geq 2$, let $c$ satisfy $|c|\leq 4$, and let $\alpha$ be a root of $cx(1-x)=1$. Then 
	$$\sum_{k=1}^\infty \frac{c^k}{k^{n}\binom{2k}{k}} \in \MZV{\{0,1,\infty,\alpha\}}{n}.$$
\end{proposition}
\begin{proof}
	The condition $|c|\leq 4$ ensures the convergence of the infinite sum. Expanding $\Li_{n-1}(x) = \sum_{k=1}^\infty x^k/k^{n-1}$ and integrating termwise, we obtain \begin{equation}\label{mzint_ex1}\int_0^1 \frac{\Li_{n-1}(cx(1-x))}{x} dx = \sum_{k=1}^\infty \frac{c^k}{k^{n}\binom{2k}{k}}. \end{equation}
	Now, $$\Li_{n-1}(cx(1-x)) = -\int_{R\circ [0,1]} \omega(0)^{n-2}\omega(1).$$
	Here $[0,1]$ denotes the path $[0,1]\to [0,1], x\mapsto x$, and $R(x)= cx(1-x)$. The above equals
	$$-\int_0^1 R^\ast \omega(0)^{n-2} R^\ast \omega(1)$$
	Now $R^\ast \omega(0) = \omega(R^{-1}(0)) - \omega(R^{-1}(\infty)) = \omega(0)+\omega(1)$ and, similarly, $R^\ast \omega(1) = \omega(\alpha) + \omega(1-\alpha)$. 
	Therefore, the series equals $$-\int_0^1 \omega(0)(\omega(0)+\omega(1))^{n-2} (\omega(\alpha)+\omega(1-\alpha)).$$
	It can be written as two iterated integrals, one supported on $\{0,1,\alpha\}$ and the other supported on $\{0,1,1-\alpha\}$. By Proposition \ref{generalStrans}, $\MZV{\{0,1,\infty,1-\alpha\}}{} = \MZV{\{0,1,\infty,\alpha\}}{}$. 
\end{proof}

The above method generalizes naturally to series that are ``twisted'' by harmonic numbers. Recall our notation $$H_{s_1,\cdots,s_k}(n) := \sum_{n>n_1,\dots,n_k\geq 1} \frac{1}{n_1^{s_1}\cdots n_k^{s_k}}.$$
\begin{proposition}\label{binomharmonictwist}
	For $n\geq 2$, let $c$ satisfy $|c|\leq 4$, and let $\alpha$ be a root of $cx(1-x)=1$. Then 
	$$\sum_{k=1}^\infty \frac{c^k H_{s_1,\cdots,s_r}(k)}{k^{n}\binom{2k}{k}} \in \MZV{\{0,1,\infty,\alpha,1-\alpha\}}{n+s_1+\cdots+s_r}.$$
\end{proposition}
\begin{proof}
	The condition $|c|\leq 4$ ensures the convergence of the infinite sum. Using the power-series expansion $\Li_{n-1,s_1,\cdots,s_r}(x) = \sum_{k=1}^\infty x^k H_{s_1,\cdots,s_r}(k)/k^{n-1}$, the infinite sum equals
	$$\int_0^1 \frac{\Li_{n-1,s_1,\cdots,s_r}(cx(1-x))}{x} dx = (-1)^{r+1}\int_0^1 \omega(0) \omega_0^{n-2}\omega_1 \omega_0^{s_1-1}\omega_1\cdots \omega_0^{s_r-1}\omega_1$$
	with $\omega_0 = \omega(0)+\omega(1),\omega_1=\omega(\alpha)+\omega(1-\alpha)$.
\end{proof}

\begin{example}
	One of the most famous special cases of the above result is the Ap{\'e}ry series
	$$\sum _{n=1}^{\infty} \frac{(-1)^n}{n^3 \binom{2 n}{n}} = -\frac{2\zeta(3)}{5}.$$
	We give another proof here. This corresponds to the case $c=-1$ and $\alpha = (1-\sqrt{5})/2$, so, by Example \ref{level5Ex}, it is a level-$5$ CMZV. Because we have found all putative $\mathbb{Q}$-relations in this space, the series evaluation follows. 
\end{example}

\begin{example}
	An equally famous example is
	$$\sum _{n=1}^{\infty} \frac{1}{n^4 \binom{2 n}{n}} = \frac{17 \pi ^4}{3240}.$$
	A mechanical proof can again be given. This corresponds to the case $c=1$ and $\alpha = e^{2\pi i /6}$, so, by Example \ref{level6Ex}, it is a level-$6$ CMZV. Since we can express level-$6$, weight-$4$ CMZVs in terms of a $\mathbb{Q}$-basis that contains $\pi^4$, this completes the proof. In fact, one need only work in $\CMZV{3}{}$ rather than $\CMZV{6}{}$; see Example \ref{level6Ex_Deligne_value}.
\end{example}

\begin{example}
	The examples above are well known, mainly because their evaluations are simple. However, our approach treats all these sums on an equal footing, regardless of the complexity of the result. We give some examples in the table below. 
	$$\begin{array}{|c|c|}
		\hline
		(c,n) & \sum _{k=1}^{\infty } \frac{c^k}{k^n \binom{2 k}{k}} \\
		\hline
		(4,3) & \pi ^2 \log (2)-\frac{7 \zeta (3)}{2} \\
		\hline
		(4,4) & 8 \text{Li}_4\left(\frac{1}{2}\right)-\frac{19 \pi ^4}{360}+\frac{\log ^4(2)}{3}+\frac{2}{3} \pi ^2 \log ^2(2) \\
		\hline
		(4,5) & -16 \text{Li}_5\left(\frac{1}{2}\right)+\pi ^2 \zeta (3)+\frac{31 \zeta (5)}{8}+\frac{2 \log ^5(2)}{15}+\frac{4}{9} \pi ^2 \log ^3(2)-\frac{19}{180} \pi ^4 \log (2) \\
		\hline
		(-1/2,3) & \frac{\log ^3(2)}{6}-\frac{\zeta (3)}{4} \\
		\hline
		(-1/2,4) & -4 \text{Li}_4\left(\frac{1}{2}\right)-\frac{13}{4} \zeta (3) \log (2)+\frac{7 \pi ^4}{180}-\frac{1}{24} 5 \log ^4(2)+\frac{1}{6} \pi ^2 \log ^2(2) \\
		\hline
		(2,3) & \pi  G-\frac{35 \zeta (3)}{16}+\frac{1}{8} \pi ^2 \log (2) \\
		\hline
		(2,4) & -2 \pi  \Im\left(\text{Li}_3\left(\frac{1}{2}+\frac{i}{2}\right)\right)+\frac{5 \text{Li}_4\left(\frac{1}{2}\right)}{2}+\frac{19 \pi ^4}{576}+\frac{5 \log ^4(2)}{48}+\frac{1}{48} \pi ^2 \log ^2(2) \\
		\hline
		(1,5) & \frac{9}{8} \pi  \sqrt{3} L_{-3}(4)+\frac{\pi ^2 \zeta (3)}{9}-\frac{19 \zeta (5)}{3} \\
		\hline
		(2-\sqrt{5},3) & 2 \text{Li}_3(\phi^{-1})-2 \zeta (3)-\frac{1}{6} \log ^3(\phi )+\frac{1}{5} \pi ^2 \log (\phi ) \\
		\hline
		(-\frac{1}{2},3) & \frac{\log ^3(2)}{6}-\frac{\zeta (3)}{4} \\
		\hline
		(-\frac{1}{2},4) & -4 \text{Li}_4\left(\frac{1}{2}\right)-\frac{13}{4} \zeta (3) \log (2)+\frac{7 \pi ^4}{180}-\frac{1}{24} 5 \log ^4(2)+\frac{1}{6} \pi ^2 \log ^2(2) \\
		\hline
		(3,3) & \frac{2 \pi  L_{-3}(2)}{\sqrt{3}}-\frac{26 \zeta (3)}{9}+\frac{2}{9} \pi ^2 \log (3) \\
		\hline
		(3,4) & -\frac{8}{3} \pi  \Im\left(\text{Li}_3\left(\frac{1}{2}+\frac{i}{2 \sqrt{3}}\right)\right)+4 \text{Li}_4\left(\frac{1}{3}\right)-\frac{\text{Li}_4\left(\frac{1}{9}\right)}{6}+\frac{29 \pi ^4}{1215}+\frac{\log ^4(3)}{18}+\frac{1}{18} \pi ^2 \log ^2(3) \\
		\hline
	\end{array}$$
	Here, $L_{-3}(s)$ is the unique primitive Dirichlet $L$-function of modulus $3$, $\phi = (\sqrt{5}+1)/2$, and $G$ is Catalan's constant. Using our Mathematica package \textsf{MultipleZetaValues}, the integral in (\ref{mzint_ex1}) can be evaluated directly with the command
	\begin{mmaCell}[functionlocal=x]{Code}
		f[c_,n_]:=\mmaDef{MZIntegrate}[PolyLog[n - 1, c*x (1 - x)]/x, {x, 0, 1}]; f[2,4]
	\end{mmaCell}

	\mycomment{For $c=4,-1/2$, the result has is CMZV of level $2$; for $c=2$, the result is CMZV of level $4$; level $5$ for $c=2-\sqrt{5}$; level $6$ for $c=-1/2,1,3$. Moreover,
		$$\sum _{n=1}^{\infty } \frac{(2-\sqrt{5})^n}{n^4 \binom{2 n}{n}} = -2 \text{Li}_3(\phi^{-1}) \log (\phi )-\frac{9}{2} \text{Li}_4(\phi^{-2}) +4 \text{Li}_4(\phi^{-1})-\frac{6}{5} \zeta (3) \log (\phi )-\frac{1}{8} 7 \log ^4(\phi )+\frac{1}{5} \pi ^2 \log ^2(\phi )$$}
\end{example}

\begin{example}
	Borwein \cite{borwein2001central} conjectured the following generalizations of Ap{\'e}ry series:
	\begin{align*}\sum _{n=1}^{\infty } \frac{(-1)^n}{n^4 \binom{2 n}{n}} &= -8 \text{Li}_3(\phi^{-1}) \log (\phi )+\frac{1}{2} \text{Li}_4(\phi^{-2})-8 \text{Li}_4(\phi^{-1}) +\frac{4}{5} \zeta (3) \log (\phi )+\frac{13 \log ^4(\phi )}{6}-\frac{7}{15} \pi ^2 \log ^2(\phi )+\frac{7 \pi ^4}{90} \\ \sum _{n=1}^{\infty } \frac{(-1)^n}{n^5 \binom{2 n}{n}} &= -\frac{5 \text{Li}_5\left(\frac{1}{\phi ^2}\right)}{2}-5 \text{Li}_4\left(\frac{1}{\phi ^2}\right) \log (\phi ) -4 \zeta (3) \log ^2(\phi )+2 \zeta (5)-\frac{4}{3}\log ^5(\phi )+\frac{4}{9} \pi ^2 \log ^3(\phi )\end{align*}
		Both sides are level-$5$ CMZVs of weights $4$ and $5$, respectively. Using our methodology, these identities can be considered established.
\end{example}

\mycomment{
	\begin{example}
		Example with many complicated $c$ are abundant, for example, let $\alpha = (\sqrt{2}+1)^4$, we have $\sum _{n=1}^{\infty } \frac{c^n}{n^s \binom{2 n}{n}} \in \CMZV{8}{s}$ with $c = \alpha^{-1}(1-\alpha)^{-1}$. When $s=4$, this is
		\small \begin{multline*}\sum _{n=1}^{\infty } \frac{1}{n^4 \binom{2 n}{n}}\left(\frac{140-99 \sqrt{2}}{8}\right)^n  = \frac{3584}{3} \sqrt{2} L(3) \log (\sqrt{2}+1)+\frac{1145 \text{Li}_4(\frac{1}{2})}{2}-\frac{1400}{3} \text{Li}_4(\frac{1}{4} (2-\sqrt{2}))\\-\frac{224 \text{Li}_4(\sqrt{2}-1)}{3}-200 \text{Li}_4(\frac{1}{\sqrt{2}})+288 \text{Li}_4(3-2 \sqrt{2})+50 \text{Li}_4(\frac{1}{8} (4-3 \sqrt{2}))-\frac{125}{12} \text{Li}_4(17-12 \sqrt{2})\\-1600 \text{Li}_3(\frac{1}{\sqrt{2}}) \log (\sqrt{2}+1)+175 \zeta (3) \log (\sqrt{2}+1)-\frac{619 \pi ^4}{432}+\frac{1205 \log ^4(2)}{192}+\frac{211}{9} \log ^4(\sqrt{2}+1)+\frac{125}{4} \log (\sqrt{2}+1) \log ^3(2)\\-50 \log ^3(\sqrt{2}+1) \log (2)-\frac{85}{24} \pi ^2 \log ^2(2)+50 \log ^2(\sqrt{2}+1) \log ^2(2)-\frac{79}{9} \pi ^2 \log ^2(\sqrt{2}+1)-\frac{200}{3} \pi ^2 \log (\sqrt{2}+1) \log (2)\end{multline*}
		here $L(s)$ is the unique primitive Dirichlet-$L$ function with modulus $8$. 
\end{example}}

\begin{example}
	When $c=1$, both $\alpha$ and $1-\alpha$ are sixth roots of unity, so every ``harmonic twist'' of $\sum \frac{1}{n^s\binom{2n}{n}}$ is a level-$6$ MZV. For example,
	$$\sum_{n=1}^\infty \frac{H_n^2}{n^2\binom{2n}{n}} = \frac{3 L_{-3}(2){}^2}{2}-\frac{\pi  L_{-3}(2) \log (3)}{\sqrt{3}}-\frac{4}{3} \pi  \Im\left(\text{Li}_3\left(\frac{1}{2}+\frac{i}{2 \sqrt{3}}\right)\right)+\frac{29 \pi ^4}{1215}+\frac{1}{36} \pi ^2 \log ^2(3)$$
	with $H_n = 1+1/2+\cdots+1/n$. 
	
	Iterated integrals supported on $\{0,1,-1,2\}$ have level $6$, so harmonic twists of $\sum \frac{(-1/2)^n}{n^s\binom{2n}{n}}$ (which has level $2$) have level $6$. For example,
	\begin{multline*}\sum _{n=1}^{\infty } \frac{(-1/2)^n H_n^2}{n^2 \binom{2 n}{n}} = \text{Li}_2\left(\frac{1}{4}\right){}^2-\frac{4 \text{Li}_4\left(\frac{1}{2}\right)}{3}+\frac{3 \text{Li}_4\left(\frac{1}{4}\right)}{2}+5 \text{Li}_2\left(\frac{1}{4}\right) \log ^2(2)-4 \text{Li}_2\left(\frac{1}{4}\right) \log (3) \log (2)+\\ 8 \text{Li}_3\left(\frac{1}{3}\right) \log (2)+2 \text{Li}_3\left(\frac{1}{4}\right) \log (2)-\frac{89}{12} \zeta (3) \log (2)-\frac{\pi ^4}{270}+\frac{77 \log ^4(2)}{18}-8 \log (3) \log ^3(2)\\ -\frac{4}{3} \log ^3(3) \log (2)+\frac{1}{18} \pi ^2 \log ^2(2) +4 \log ^2(3) \log ^2(2)+\frac{2}{3} \pi ^2 \log (3) \log (2)\end{multline*}
	
	Harmonic twists of the Ap{\'e}ry series $\sum \frac{(-1)^n}{n^s\binom{2n}{n}}$ are level-$10$ CMZVs; see Example \ref{level10Ex1}. We give an example whose simplicity remains unexplained:
	$$\sum _{n=1}^{\infty } \frac{(-1)^n H_n}{n^3 \binom{2 n}{n}} = -\frac{12}{5} \text{Li}_3(\phi^{-1}) \log (\phi )+\frac{3}{20} \text{Li}_4(\phi^{-2}) -\frac{12}{5} \text{Li}_4(\phi^{-1})+\frac{6}{25} \zeta (3) \log (\phi )+\frac{13 \log ^4(\phi )}{20}-\frac{7}{50} \pi ^2 \log ^2(\phi )+\frac{\pi ^4}{50}.$$
\end{example}
One could give many more examples, but we stop here. \par

The statement and proof of Theorem \ref{binomharmonictwist} cover only the case $n\geq 2$. What happens when $n=1$? An analogue holds after tensoring with the field $\mathbb{Q}(\alpha)$. 
\begin{theorem}
	Let $c$ satisfy $|c|\leq 4$, and let $\alpha$ be a root of $cx(1-x)=1$. Then 
	$$\sum_{k=1}^\infty \frac{c^k H_{s_1,\cdots,s_r}(k)}{k\binom{2k}{k}} \in \MZV{\{0,1,\infty,\alpha,1-\alpha\}}{1+s_1+\cdots+s_r}\otimes_\mathbb{Q} \mathbb{Q}(\alpha).$$
\end{theorem}
Although the proof is not difficult, we defer it, along with further examples, to a later article.

\begin{example}
	Sun \cite{sun2010conjectures} conjectured that \begin{equation}\label{sun1}\sum_{n=1}^\infty \frac{1}{n^2\binom{2n}{n}} (3H_{n-1}^2+\frac{4}{n}H_{n-1}) = \frac{\pi^4}{360}.\end{equation}
	By Proposition \ref{binomharmonictwist}, the left-hand side lies in $\MZV{\{0,1,\infty,\alpha,1-\alpha\}}{4}$, where $\alpha$ is a solution of $\alpha(1-\alpha)=1$, and thus lies in $\CMZV{6}{4}$. Hence, even without any calculation, we know that the above conjecture can be proved using our machinery. We give more details for this illustrative example. One observes that $$\sum_{n=1}^\infty (3H_{n-1}^2 + \frac{4}{n} H_{n-1})z^n = 4 \Li_{2,1}(z)+3\Li_{1,2}(z)+6\Li_{1,1,1}(z).$$
	Here, we have used our notation for generalized polylogarithms.
	Hence, the sum equals $$\int_0^1 \frac{4 \Li_{2,1}(x(1-x))+3\Li_{1,2}(x(1-x))+6\Li_{1,1,1}(x(1-x))}{x} dx.$$
	Recall that $\int_0^1 x^{n-1}(1-x)^n = \frac{1}{n}\binom{2n}{n}^{-1}$.
	Each integral can be written as an element of $\MZV{\{0,1,\infty,\alpha,1-\alpha\}}{4}$. For example, the first term equals $\int_0^1 4\omega(0)\omega_0 \omega_1\omega_1$, with $\omega_0=\omega(0)+\omega(1)$ and $\omega_1 = \omega(\alpha)+\omega(1-\alpha)$, where $\alpha$ is a sixth root of unity. Converting all the above terms into level-$6$, weight-$4$ CMZVs proves the result.\par
	The above procedure is executed automatically by our Mathematica package with the following command:
	\begin{mmaCell}[functionlocal=x]{Code}
		\mmaDef{MZIntegrate}[(4 \mmaDef{MZPolyLog}[{0, 1, 1}, x (1 - x)] + 3 \mmaDef{MZPolyLog}[{1, 0, 1}, x (1 - x)] 
		+6 \mmaDef{MZPolyLog}[{1, 1, 1}, x (1 - x)])/x, {x, 0, 1}]
	\end{mmaCell}
	\begin{mmaCell}{Output}
		Pi^4/360
	\end{mmaCell}
\end{example}

We mentioned another integration kernel:
$$\int_0^1 \frac{x^n(1-x)^n}{x}(-\log x) dx = \frac{1}{n\binom{2n}{n}} (H_{2n}-H_{n-1}),$$
which follows because the left-hand side is a derivative of Euler's beta function. This enables us, for example, to express the series $$\sum _{n=1}^{\infty } \frac{2^n \left(-3 H_n+2 H_{2 n}+\frac{2}{n}\right)}{n^2 \binom{2 n}{n}}$$
as $$\int_0^1 \frac{-\Li_{1,1}(2 x (1-x))-\text{Li}_2(2 x (1-x))+2 \text{Li}_1(2 x (1-x)) (-\log (x))}{x} dx.$$
The roots of $2x(1-x)=1$ are $(1\pm i)/2$, so this expression is in $\MZV{S}{3}$ for $S=\{0,1,\infty,(1+i)/2,(1-i)/2\}$, which is contained in $\CMZV{4}{3}$. We have therefore proved the first of the following identities.

\begin{proposition}[Conjectures from \cite{sun2010conjectures}]
	The following are all true:
	
	\noindent\begin{minipage}{.5\linewidth}
		$$\begin{aligned}
			\sum _{n=1}^{\infty } \frac{2^n \left(-3 H_n+2 H_{2 n}+\frac{2}{n}\right)}{n^2 \binom{2 n}{n}} &=\frac{7 \zeta (3)}{4},\\
			\sum _{n=1}^{\infty } \frac{2^n \left(-11 H_n+6 H_{2 n}+\frac{8}{n}\right)}{n^2 \binom{2 n}{n}}&=2 \pi  G,\\
			\sum _{n=1}^{\infty } \frac{3^n \left(-10 H_n+6 H_{2 n}+\frac{7}{n}\right)}{n^2 \binom{2 n}{n}}&=2 \sqrt{3} \pi  L_{-3}(2),\\
			\sum _{n=1}^{\infty } \frac{3^n \left(H_n+\frac{1}{2 n}\right)}{n^2 \binom{2 n}{n}}&=\frac{1}{3} \pi ^2 \log (3), \\
			\sum _{n=1}^{\infty }\frac{17 H_n+H_{2 n}}{n^2 \binom{2 n}{n}}&=\frac{5}{2} \sqrt{3} \pi  L_{-3}(2), \\
			\sum _{n=1}^{\infty } \frac{H_3^{(n)}}{n^2 \binom{2 n}{n}}&=\frac{\pi ^2 \zeta (3)}{27}+\frac{\zeta (5)}{9}, \\
			\sum _{n=1}^{\infty } \frac{-163 H_n+97 H_{2 n}+\frac{227}{n}}{n^4 \binom{2 n}{n}}&=\frac{165}{8} \sqrt{3} \pi  L_{-3}(4),
		\end{aligned}$$
	\end{minipage}%
	\begin{minipage}{.5\linewidth}
		$$\begin{aligned}
			\sum _{n=1}^{\infty } \frac{2^n \left(-7 H_n+2 H_{2 n}+\frac{2}{n}\right)}{n^2 \binom{2 n}{n}}&=-\frac{\pi^2}{2} \log (2),\\
			\sum _{n=1}^{\infty } \frac{3^n \left(-8 H_n+6 H_{2 n}+\frac{5}{n}\right)}{n^2 \binom{2 n}{n}}&=\frac{26 \zeta (3)}{3},\\
			\sum _{n=1}^{\infty } \frac{H_{2 n}+\frac{2}{3 n}}{n^2 \binom{2 n}{n}}&=\zeta (3),\\
			\sum _{n=1}^{\infty } \frac{2 H_n+H_{2 n}}{n^2 \binom{2 n}{n}}&=\frac{5 \zeta (3)}{3},\\
			\sum _{n=1}^{\infty } \frac{-H_n+H_{2 n}+\frac{2}{n}}{n^4 \binom{2 n}{n}}&=\frac{11 \zeta (5)}{9},\\
			\sum _{n=1}^{\infty } \frac{-102 H_n+3 H_{2 n}+\frac{28}{n}}{n^4 \binom{2 n}{n}}&=-\frac{55}{18}  \pi ^2 \zeta (3).
		\end{aligned}$$
	\end{minipage}
	Here, $G$ is Catalan's constant, and $L_{-3}(s)$ is the unique primitive Dirichlet $L$-function of modulus $3$.
\end{proposition}
\begin{proof}
	Using the method above, it is straightforward to convert these sums into integrals. For example, the first, second, third, and penultimate identities correspond to the integrals
	$$\begin{aligned}&\int_0^1 \frac{\Li_1(2x(1-x))(-2\log x) - \Li_2(2x(1-x))-\Li_{1,1}(2x(1-x))}{x} dx \\ 
		&\int_0^1 \frac{\Li_1(2x(1-x))(-6\log x) - 3\Li_2(2x(1-x))-5\Li_{1,1}(2x(1-x))}{x} dx \\
		&\int_0^1 \frac{\Li_1(2x(1-x))(-2\log x) - 5\Li_2(2x(1-x))-5\Li_{1,1}(2x(1-x))}{x} dx \\
		&\int_0^1 \frac{\Li_3(x(1-x))(-3\log x) - 99\Li_{3,1}(x(1-x))-74\Li_{4}(x(1-x))}{x} dx
	\end{aligned}$$
	They can all be converted into CMZVs of level $4$ or $6$ and weight $\leq 5$, thus giving the claimed evaluations. All of these identities were already established by Ablinger \cite{ablinger2014iterated,ablinger2011harmonic,ablinger2019proving,ablinger2017discovering}.
\end{proof}

\begin{proposition}[Conjectures in \cite{sun2010conjectures}]
	The following identities hold:
	$$\begin{aligned}\sum _{n=1}^{\infty } \frac{\left(\left(\frac{\sqrt{5}-1}{2} \right)^{2 n}+\left(\frac{\sqrt{5}+1}{2}\right)^{2 n}\right) \left(H_{2 n}-H_{n-1}\right)}{n^2 \binom{2 n}{n}} &=\frac{4}{25} \pi ^2 \log  \phi +\frac{41 \zeta (3)}{25} \\
		\sum _{n=1}^{\infty } \frac{\left(\left(\frac{5-\sqrt{5}}{2}\right)^{n}+\left(\frac{5+\sqrt{5}}{2} \right)^{n}\right) \left(H_{2 n}-H_{n-1}\right)}{n^2 \binom{2 n}{n}}&=\frac{62 \zeta (3)}{25}+\frac{3}{25} \pi ^2 \log \phi+\frac{1}{10} \pi ^2 \log 5
	\end{aligned}$$
	Here, $\phi = (\sqrt{5}+1)/2$.
\end{proposition}
\begin{proof}
	For the first identity, let $\alpha$ be a root of $((\sqrt{5}-1)/2)^2 x(1-x)=1$; for the second, let it be a root of $(5-\sqrt{5})x(1-x)/2=1$. In either case, set $S=\{0,1,\infty,\alpha,1-\alpha\}$. Then $\MZV{S}{}\subset \CMZV{5}{}$.
\end{proof}

The above two identities have been proved in \cite{xu2022sun,xu2022note}. Our approach can also evaluate
$$\sum _{n=1}^{\infty } \frac{\left(\frac{\sqrt{5}+1}{2}\right)^{2 n} \left(H_{2 n}-H_{n-1}\right)}{n^2 \binom{2 n}{n}} = -\int_0^1 \frac{\text{Li}_1\left((\frac{\sqrt{5}+1}{2} )^2 x (1-x)\right) \log (x)}{x} dx.$$
The result is
$$\frac{3}{5} \text{Li}_3\left(\frac{\sqrt{5}-1}{2}\right)+\frac{19 \zeta (3)}{25}-\frac{1}{5} \log ^3(\phi )+\frac{6}{25} \pi ^2 \log (\phi ).$$

The next example uses the integration kernel $$\int_0^1 \frac{(x^2(1-x))^n}{x}\log(\frac{x}{1-x}) dx = \frac{1}{2n\binom{3n}{n}} (H_{2n-1}-H_n).$$
\begin{proposition}[Conjecture 10.61 in \cite{sun2021book}]
	For $s=1,2$, respectively, $$\sum _{n=1}^{\infty } \frac{H_{2 n}-H_n}{\binom{3 n}{n} \left(2^n n^s\right)}$$ equals $$-\frac{\pi ^2}{60}+\frac{3 \log ^2(2)}{10}+\frac{1}{20} \pi  \log (2) \qquad \text{and} \qquad -\frac{\pi  G}{2}+\frac{33 \zeta (3)}{32}+\frac{1}{24} \pi ^2 \log (2).$$
\end{proposition}
\begin{proof}
	We have $$\sum _{n=1}^{\infty } \frac{H_{2 n}-H_n}{\binom{3 n}{n} 2^n n^s}=\int_0^1 \frac{\text{Li}_s\left(\frac{1}{2} x^2 (1-x)\right)+2 (\log (x)-\log (1-x)) \text{Li}_{s-1}\left(\frac{1}{2} x^2 (1-x)\right)}{x} \, dx$$
	Here, we interpret $\Li_0(x) = x/(1-x)$. Since $\frac{x^2(1-x)}{2}=1$ implies that $x=-1,1\pm i$, the right-hand side can be written as an iterated integral in which $\omega(-1)$, $\omega(1-i)$, and $\omega(1+i)$ do not occur in the monomial. Moreover, $\MZV{S_1}{s+1},\MZV{S_2}{s+1},\MZV{S_3}{s+1} \subset \CMZV{4}{s+1}$, where $$S_1 = \{0,1,\infty,-1\}, \qquad S_2=\{0,1,\infty,1-i\}, \qquad S_3=\{0,1,\infty,1+i\},$$
	completing the proof. The procedure above is automatically performed with the Mathematica command 
	
	\begin{mmaCell}[functionlocal=x]{Code}
		\mmaDef{MZIntegrate}[(2PolyLog[#-1,1/2x^2(1-x)](Log[x]-Log[1-x])
		+PolyLog[#,1/2x^2(1-x)])/x,{x,0,1}]&/@{1,2}
	\end{mmaCell}
\end{proof}

Nothing prevents us from taking $s=3$ or larger. For example, we have
$$\sum _{n=1}^{\infty } \frac{H_{2 n}-H_n}{\binom{3 n}{n} 2^n n^3} = \frac{1}{2} \pi  G \log (2)+\frac{9 \text{Li}_4\left(\frac{1}{2}\right)}{2}+\frac{93}{32} \zeta (3) \log (2)-\frac{31 \pi ^4}{640}+\frac{3 \log ^4(2)}{16}-\frac{5}{24} \pi ^2 \log ^2(2)$$
The MZV nature of the simpler series $\sum \frac{1}{\binom{3n}{n}2^n n^s}$ was already revealed in the author's previous paper \cite{au2020evaluation}. Borwein \cite{borwein2004experimentation} investigated these series experimentally. 

\mycomment{
	\begin{proposition}[Conjectures in \cite{sun2010conjectures}]The followings are true:
		$$\begin{aligned}\sum _{n=1}^{\infty } \frac{(-1)^{n-1} \left(10 H_n-\frac{3}{n}\right)}{n^3 \binom{2 n}{n}} &=\frac{\pi ^4}{30} \\
			\sum _{n=1}^{\infty } \frac{(-1)^{n-1} \left(4 H_n+H_{2 n}\right)}{n^3 \binom{2 n}{n}} &=\frac{2 \pi ^4}{75}
		\end{aligned}$$
	\end{proposition}
	\begin{proof}
		Using the same technique as above, one can show above two series are in $\MZV{S}{4}$ with $S=\{(1-\sqrt{5})/2,(1+\sqrt{5})/2,0,1,\infty\} \subset \CMZV{10}{4}$. CMZV of level 10 weight 4 has not been solved\footnote{We are 2 relations away from Deligne's bound (Theorem \ref{deligne_bound}), see remarks at end of \ref{nonsd}}, nonetheless, the currently available relations suffice to deduce them.
	\end{proof}
}

We give an example that involves a multiple polylogarithm rather than a generalized polylogarithm. 
\begin{proposition}[Conjecture in \cite{sun2021book}]\label{multipolylogconj}
	Let $\lfloor x \rfloor$ denote the floor function. Then
	\begin{equation}\sum _{n=1}^{\infty } \frac{2^n \left(H_{\left\lfloor n/2\right\rfloor }-\frac{2 (-1)^n}{n}\right)}{n^2 \binom{2 n}{n}}=\frac{7 \zeta (3)}{4}.\end{equation}
\end{proposition}
\begin{proof}
	Note that $2x(1-x)=1$ if and only if $x=\frac{1\pm i}{2}$, while $-2x(1-x) = 1$ if and only if $x=(1\pm \sqrt{3})/2$. Let $$S = \{0,1,\infty,\frac{1\pm i}{2},(1\pm \sqrt{3})/2\}.$$ It can be shown, as in Corollary \ref{main_corollary}, that $\MZV{S}{} \subset \CMZV{12}{}$. Since
	$$H_{\left\lfloor n/2\right\rfloor } = H_{n-1}+a_{n-1} + \frac{1+(-1)^n}{2}$$
	where $a_n = \sum_{k=1}^n (-1)^k/k$. The infinite series for the first and last terms are in $\MZV{S}{3}$. It remains to handle $A = \sum_{n=1}^\infty \frac{2^n}{n^2\binom{2n}{n}} a_{n-1}$. First, note that $a_{n-1}$ is the coefficient in the multiple polylogarithm $$\Li_{1,1}(x,-1) = \sum_{n\geq 1} \frac{x^n a_{n-1}}{n} = \int_0^x \omega(1)\omega(-1).$$
	Therefore, by the pullback formula for iterated integrals, $$\Li_{1,1}(2x(1-x),-1) = \int_0^x R^\ast \omega(1) R^\ast \omega(-1) = \int_0^x (\omega(\frac{1+i}{2})+\omega(\frac{1-i}{2}))(\omega(\frac{1+\sqrt{3}}{2})+\omega(\frac{1-\sqrt{3}}{2})),$$
	where $R(x) = 2x(1-x)$.
	Hence, $$A = \int_0^1 \frac{\Li_{1,1}(2x(1-x),-1)}{x} dx = \int_0^1 \omega(0)(\omega(\frac{1+i}{2})+\omega(\frac{1-i}{2}))(\omega(\frac{1+\sqrt{3}}{2})+\omega(\frac{1-\sqrt{3}}{2})).$$
	Therefore, $A\in \CMZV{12}{3}$. Since level-$12$, weight-$3$ CMZVs are included in the datamine, the claim follows.\footnote{The integral $A$ can be calculated with the Mathematica package using the command \texttt{IterInt}.} 
\end{proof}

\bibliographystyle{plain} % We choose the "plain" reference style
\bibliography{polylog1_2021} % Paper-specific, audited bibliography

\end{document}